\documentclass[a4paper,10pt,oneside]{amsart}

\RequirePackage{mathmacros}
\RequirePackage{upref}
\RequirePackage{amsthm}
\RequirePackage{enumerate}
\usepackage{xr-hyper}
\usepackage[hyperindex=true, pdfauthor= Claus\  Gerhardt, pdftitle= Scalar\ Curvature, bookmarks=true, extension= pdf, colorlinks=true, plainpages=false,hyperfootnotes=false, debug=false, pagebackref]{hyperref}

\theoremstyle{plain}

\theoremstyle{definition}

\swapnumbers
\theoremstyle{remark}

\numberwithin{equation}{section}

\begin{document}
\title[Prescribed Scalar Curvature]{Hypersurfaces of prescribed Scalar curvature in
Lorentzian manifolds}

\author{Claus Gerhardt}
\address{Ruprecht-Karls-Universit\"at, Institut f\"ur Angewandte Mathematik,
Im Neuenheimer Feld 294, 69120 Heidelberg, Germany}
\email{gerhardt@math.uni-heidelberg.de}
\urladdr{\url{http://www.math.uni-heidelberg.de/studinfo/gerhardt/}}

%
\date{September 30, 2001}
%

\dedicatory{Dedicated to Robert Finn on the occasion of his eightieth birthday}

\begin{abstract} The existence of closed hypersurfaces of prescribed sca\-lar
curvature  in globally hyperbolic Lorentzian manifolds is proved provided there
are barriers.
\end{abstract}
\maketitle
\thispagestyle{empty}

\tableofcontents
\setcounter{section}{-1}
\section{Introduction} 

Consider the problem of finding a closed hypersurface of prescribed curvature
$F$ in a globally hyperbolic \di {(n+1)} Lorentzian mani\-fold $N$ having a
compact Cauchy hypersurface $\so$. To be more precise, let
$\Om$ be a connected open subset of $N,f\in C^{2,\al}(\bar\Om), F$ a smooth,
symmetric function defined in an open cone $\C\su\R[n]$, then we look for a
space-like hypersurface
$M\su \Om$ such that
\begin{equation}\lae{0.1}
F_{|_M}=f(x)\qq\A x\in M,
\end{equation}
where $F_{|_M}$ means that $F$ is evaluated at the vector $(\ka_i(x))$ the
components of which are the principal curvatures of $M$. The prescribed function
$f$ should satisfy natural structural conditions, e. g. if $\varGamma$ is the
positive cone and the hypersurface $M$ is supposed to be convex, then $f$ should
be positive, but no further, merely technical, conditions should be imposed.

In \ci{rb, rb2, cg1, cg7}   the case $F=H$, the mean curvature, has been treated, 
and in
\ci{cg8} we solved the problem for curvature functions $F$ of class $K^*$ that
includes the Gaussian curvature, see \ci[Section 1]{cg8} for the definition, but
excludes the symmetric polynomials $H_k$ for $1<k<n$. Among these, $H_2$,
that corresponds to the scalar curvature operator, is of special interest.

However, a solution of equation \re{0.1} with $F=H_2$ is in general not a
hypersurface of prescribed scalar curvature---unless the ambient space has
constant curvature---since the scalar curvature of a hypersurface also
depends on $\bar R_{\al\bet}\n^a\n^\bet$. Thus, we have to allow that the
right-hand side $f$ also  depends on  time-like vectors  and look for
hypersurfaces $M$ satisfying

\begin{equation}\lae{0.2}
F_{|_M}=f(x,\n)\qq\A x\in M,
\end{equation}

\cvm
\nd where $\n=\n(x)$ is the past-directed normal of $M$ in the point $x$.

\cvb
To give a precise statement of the existence result we need a few definitions and
assumptions. First, we assume that $\Om$ is a precompact, connected, 
open subset of $N$, that is bounded by  two {\it achronal}, connected, space-like
hypersurfaces $M_1$ and $M_2$ of class $C^{4,\al}$, where $M_1$ is supposed
to lie in the past of
$M_2$.

\cvb
Let $F=H_2$ be the scalar curvature operator defined on the open cone
$\C_2\su \R[n]$, and
$f=f(x,\n)$ be of class
$C^{2,\al}$ in its arguments such that 

\begin{align}
0<c_1&\le f(x,\n)\qq\tup{if}\q\spd\n\n=-1,\lae{0.3}\\[\cma]
\nnorm{f_\bet(x,\n)}&\le c_2 (1+\nnorm\n^2),\lae{0.4}\\
\intertext{and}
\nnorm{f_{\n^\bet}(x,\n)}&\le c_3 (1+\nnorm\n),\lae{0.5}
\end{align}

\cvm
\nd for all $x\in\bar\Om$ and all past directed time-like vectors $\n\in T_x(\Om)$,
where $\nnorm{\cdot}$ is a Riemannian reference metric that will be detailed in
\rs{2}. 

\cvm
We suppose that the boundary components $M_i$ act as barriers for  $(F,f)$.

\cvb
\bd
$M_2$ is an {\it upper barrier} for $(F,f)$, if $M_2$ is {\it admissible}, i.e. its
principal curvatures $(\ka_i)$ with respect to the past directed normal belong to
$\C_2$, and if

\begin{equation}
\fv F{M_2}\ge f(x,\n)\qq\A\,x\in M_2.
\end{equation}

\cvm
$M_1$ is a lower barrier for $(F,f)$, if at the points $\Si\su M_1$, where $M_1$
is admissible, there holds

\begin{equation}
\fv F\Si \le f(x,\n)\qq\A\,x\in \Si.
\end{equation}

\cvm
\nd $\Si$ may be empty.
\ed

\cvb
\br
This definition of upper and lower barriers for a pair $(F,f)$ also makes sense for
other curvature functions $F$ defined in an open convex cone $\C$, with a
corresponding meaning of the notion {\it admissable}.
\er

\cvb
Now, we can state the main theorem.

\cvb
\bt\lat{0.2}
Let $M_1$ be a lower and $M_2$ an upper barrier for $(F,f)$, where $F=H_2$. Then, the problem

\begin{equation}
\fmo M= f(x,\n)
\end{equation}

\cvm
\nd
has an admissible solution $M\su \bar\Om$ of class $C^{4,\al}$ that can be
written as a graph over $\mc S_0$ provided there exists a strictly convex
function $\chi\in C^2(\bar\Om)$.
\et

\cvb
\br
As we have shown in \ci[Lemma 2.7]{cg8} the existence of a strictly convex
function
$\chi$ is guaranteed by the assumption that the level hypersurfaces
$\{x^0=\tup{const}\}$ are strictly convex in $\bar\Om$, where $(x^\al)$ is a
Gaussian coordinate system associated with $\so$.

Looking at  Robertson-Walker space-times it seems that the assumption of the
existence of  a strictly convex function in the neighbourhood of a given compact
set is not too restrictive: in Minkowski space e.g. $\chi=-\abs{x^0}^2 +\abs x^2$
is a globally defined strictly convex function. The only obstruction we are aware
of is the existence of a compact maximal slice. In the neighbourhood of such  a
slice a strictly convex function cannot exist.
\er

\cvb
The existence result of our main theorem would also be valid in Riemannian
manifolds if one could prove $C^1$- estimates. For the $C^2$- estimates the
nature of the ambient space is irrelevant though the proofs are slightly
different.

\cvm
For prescribed curvature problems it seems more natural to assume that the
right-hand side $f$ depends on $(x,\n)$, and we shall prove in a subsequent
paper existence results for curvature functions $F\in (K^*)$, where the ambient
space can be Riemannian or Lorentzian, cf. \ci{cg9}.

\cvb
The paper is organized as follows: In \rs 1 we take a closer look at   curvature
functions and define the concept of {\it elliptic regularization} for these functions,
and analyze some of its properties.

\cvm

In \rs 2 we introduce the notations and common definitions we rely on,  and state
the equations of Gau{\ss}, Codazzi, and Weingarten for space-like hypersurfaces.

\cvm
In \rs 3 we look at the curvature flow associated with our problem, and the
corresponding evolution equations for the basic geometrical quantities of the
flow hypersurfaces.

\cvm
In \rs 4 we prove lower order estimates for the evolution problem,  while a
priori estimates in the
$C^2$\nbd{norm} are derived in
\rs 5.

\cvm
In \rs 6, we demonstrate that the evolutionary solution converges to a
stationary approximation of our problem, i.e. to a solution for a curvature
problem, where
$F$ is replaced by its elliptic regularization $F_\e$.

\cvm
The uniform $C^1$- estimates for the stationary approximations are derived in
Sections \ref{S:7} and \ref{S:8}, the $C^2$- estimates are given in \rs{9}, while
the final existence result is contained in \rs{10}.

\cvb
\section{Curvature functions}\las{1}

\cvb
Let $\C\su\R[n]$ be an open  cone containing the positive cone $\C_+$, and
$F\in C^{2,\al}(\C)\ii C^0(\bar\C)$ a positive symmetric function satisfying the
condition
\begin{equation}\lae{1.1}
F_i=\pd F\ka i>0\; ;
\end{equation}
then, $F$ can also be viewed as a function defined on the space of symmetric
 matrices $\mathscr C$, the eigenvalues of which belong to $\C$, namely, let
$(h_{ij})\in
\msc C$ with eigenvalues $\ka_i,\,1\le i\le n$, then define $F$ on $\msc C$
by
\begin{equation}
F(h_{ij})=F(\ka_i).
\end{equation}

If we define 
\begin{align}
F^{ij}&=\pde F{h_{ij}}\\
\intertext{and}
F^{ij,kl}&=\pddc Fh{{ij}}{{kl}}
\end{align}
then, 
\begin{equation}
F^{ij}\x_i\x_j=\pdc F\ka i \abs{\x^i}^2\q\A\, \x\in\R[n],
\end{equation}
\nd in an appropriate coordinate system,

\begin{equation}
F^{ij} \,\text{is diagonal if $h_{ij}$ is diagonal,}
\end{equation}
and
\begin{equation}\lae{1.7}
F^{ij,kl}\h_{ij}\h_{kl}=\pddc F{\ka}ij\h_{ii}\h_{jj}+\sum_{i\ne
j}\frac{F_i-F_j}{\ka_i-\ka_j}(\h_{ij})^2,
\end{equation}
for any $(\h_{ij})\in \msc S$, where $\msc S$ is the space of all symmetric
matrices. The second term on the right-hand side of \re{1.7} is non-positive if
$F$ is concave, and non-negative if $F$ is convex, and has to be interpreted as a
limit if $\ka_i=\ka_j$.

The preceding considerations are
also applicable if the
$\ka_i$ are the principal curvatures of a space-like hypersurface $M$ with metric
$(g_{ij})$.
$F$ can then be looked at as being defined on the space of all symmetric tensors
$(h_{ij})$ the eigenvalues of which 
belong to
$\C$. Such tensors will be called \tit{admissible}; when the second
fundamental form of $M$ is admissible, then, we also call $M$ admissible.

\cvb
For an admissible tensor $(h_{ij})$

\begin{equation}
F^{ij}=\pdc Fh{{ij}}
\end{equation}

\cvm
\nd
is  a contravariant tensor of second order. Sometimes it will be convenient
to circumvent the dependence on the metric by considering $F$ to depend on the
mixed tensor

\begin{equation}
h_j^i=g^{ik}h_{kj}.
\end{equation}

Then,

\begin{equation}
F_i^j=\pdm Fhji
\end{equation}

\cvm
\nd
is also a mixed tensor with contravariant index $j$ and covariant index $i$.

\cvm
Such functions $F$ are called curvature functions. Important examples are the
symmetric polynomials of order $k$, $H_k$, $1\le k\le n$,

\begin{equation}
H_k(\ka_i)=\sum_{i_1<\cdots <i_k} \ka_{i_1}\cdots \ka_{i_k}.
\end{equation}

They are defined on an open cone $\C_k$  that can be
characterized as the connected component of $\{H_k>0\}$ that contains
$\C_+$.

\cvm
Since we have in mind that the $\ka_i$ are the principal curvatures of a
hypersurface, we use the standard symbols $H$ and $\abs A$ for

\begin{align}
H&=\sum_i\ka_i,\\
\intertext{and}
\abs A^2&=\sum_i \ka^2_i.
\end{align}

\cvm
The scalar curvature function $F=H_2$ can then be expressed as

\begin{equation}
F=\frac12(H^2-\abs A^2),
\end{equation}

\cvm
\nd and we deduce that for $(\ka_i)\in \C_2$

\begin{align}
\abs A^2&\le H^2,\lae{1.15}\\[\cma]
F_i&=H-\ka_i,\lae{1.16}\\
\intertext{and hence,}
H F_i&\ge F,\lae{1.17}
\end{align}

\cvm
\nd for \re{1.17} is equivalent to

\begin{equation}
H\ka_i\le \frac12 H^2 + \frac12 \abs A^2,
\end{equation}

\cvm
\nd which is obviously valid.

\cvb
In important ingredient in our existence proof will be the method of {\it elliptic
regularization}

\cvm
\bl\lal{1.1}
For each $\e>0$, consider the linear isomorphism $\f_\e$ in $\R[n]$ given by

\begin{equation}\lae{1.21}
( \tilde \ka_i)=\f_\e(\ka_i)=(\ka_i+\e H).
\end{equation}
Let $F\in C^2(\C)\ii C^0(\bar\C)$ be a curvature function such that

\begin{equation}\lae{1.19}
\fv F{\pa\C}=0.
\end{equation}

\cvm\nd Then, $\C_\e=\f_\e^{-1}(\C)$ is an open cone and
$F_\e=F\circ\f_\e\in C^2(\C_\e)\ii C^0(\bar\C_\e)$ a curvature function
satisfying

\begin{equation}
\fv {F_\e}{\pa\C_\e}=0.
\end{equation}

\cvm
\nd Assume furthermore, that

\begin{equation}\lae{1.20}
H>0\qq\tup{in}\q\C.
\end{equation}

\cvm
\nd Then,

\begin{align}\lae{1.23b}
\C&\su \C_\e,\\
\intertext{and}
H&>0\qq\tup{in}\q \C_\e.\lae{1.24b}
\end{align}
\el

\bp
We only prove the assertions \re{1.23b} and \re{1.24b} since the other
assertions are  obvious. Let $(\ka_i)\in\C$ be fixed. Then,

\begin{equation}
0< F(\ka_i)\le F(\ka_i+\e H),
\end{equation}

\cvm
\nd because $F$ is monotone, and we deduce

\begin{equation}
(\ka_i+\e H)\in \C\qq\A\,\e>0,
\end{equation}

\cvm
\nd in view of  \re{1.20} and the monotonicity of $F$, cf. \re{1.1}.

\cvm
To prove \re{1.24b}, we observe that

\begin{equation}
\sum_i \tilde \ka_i=(1+\e n)\sum_i\ka_i.
\end{equation}
\ep

\cvb
\br
(i) Let $F$ be as in \rl{1.1} and assume moreover, that $F$ is homogeneous of
degree $1$, and concave, then, 

\begin{equation}\lae{1.25}
F(\ka_i)\le \frac1n F(1,\dotsc, 1) H\qq\A\,(\ka_i)\in\C,
\end{equation}

\cvm
\nd and we conclude that condition \re{1.20} is satisfied.

\cvm
(ii) Let $F$ be as in \rl{1.1}, but suppose that $F$ is homogeneous of degree
$d_0>0$ and $F^\frac1{d_0}$ concave, then, the relation \re{1.20} is also valid.
\er

\cvb
\bp
The inequality \re{1.25} follows easily from the concavity and homogeneity

\begin{equation}
\begin{aligned}
F(\ka_i)&\le  F(1,\dotsc, 1)+ \sum_iF_i(1,\dotsc, 1) (\ka_i-1)\\[\cma]
&= \frac1n
F(1,\dotsc, 1)H,
\end{aligned}
\end{equation}

\cvm
\nd since $F_i(1,\dotsc, 1)=\frac1n F(1,\dotsc, 1)$, while the other assertions
are obvious.
\ep

\cvb
For better reference, we use a tensor setting in the next lemma, i.e. the
$(\ka_i)\in\C$ are the eigenvalues of an admissible tensor $(h_{ij})$ with respect
to a Riemannian metric $(g_{ij})$. In this setting the elliptic regularization of $F$
is given by

\begin{equation}
\tilde F(h_{ij})\equiv F(h_{ij}+\e H g_{ij}).
\end{equation}

\cvb
\bl\lal{1.4}
Let $ \tilde F$ be the elliptic regularization of a curvature function $F$ of class
$C^2$, then,

\begin{equation}\lae{1.28}
\tilde F^{ij}=F^{ij}+\e F^{rs}g_{rs} g^{ij},
\end{equation}

\cvm
\nd and

\begin{equation}\lae{1.29}
\begin{aligned}
\tilde F^{ij,kl}&= F^{ij,kl} +\e F^{ij,ab} g_{ab} g^{kl}\\[\cma]
&\hp{=\; }+ \e F^{rs,kl} g_{rs} g^{ij} +\e^2 F^{rs,ab} g_{rs} g_{ab} g^{ij} g^{kl}.
\end{aligned}
\end{equation}

\cvm
\nd If $F$ is concave, then, $ \tilde F$ is also concave.
\el

\cvb
\bp
The relations \re{1.28} and \re{1.29} are straight-forward calculations.

\cvm
To prove the concavity of $ \tilde F$, let $(\h_{ij})$ be a symmetric tensor,
then,

\begin{equation}
\begin{aligned}
\tilde F^{ij,kl} \h_{ij} \h_{kl}&= F^{ij,kl} \h_{ij}\h_{kl} +2\e F^{ij,rs} \h_{ij}
g_{rs} g^{kl}\h_{kl}\\[\cma]
&\hp{=\; }+ \e^2 F^{rs,ab} g_{rs} g_{ab} (g^{ij}\h_{ij})^2\le 0.
\end{aligned}
\end{equation}
\ep

\cvb
\section{Notations and preliminary results}\las{2}

\cvb
The main objective of this section is to state the equations of Gau{\ss}, Codazzi,
and Weingarten for space-like hypersurfaces $M$ in a  \di{(n+1)} Lorentzian 
space
$N$. Geometric quantities in $N$ will be
denoted by
$(\bar g_{\al\bet}),(\riema \al\bet\ga\de)$, etc., and those in $M$ by $(g_{ij}), (\riem
ijkl)$, etc. Greek indices range from $0$ to $n$ and Latin from $1$ to $n$; the
summation convention is always used. Generic coordinate systems in $N$ resp.
$M$ will be denoted by $(x^\al)$ resp. $(\x^i)$. Covariant differentiation will
simply be indicated by indices, only in case of possible ambiguity they will be
preceded by a semicolon, i.e. for a function $u$ in $N$, $(u_\al)$ will be the
gradient and
$(u_{\al\bet})$ the Hessian, but e.g., the covariant derivative of the curvature
tensor will be abbreviated by $\riema \al\bet\ga{\de;\e}$. We also point out that

\begin{equation}
\riema \al\bet\ga{\de;i}=\riema \al\bet\ga{\de;\e}x_i^\e
\end{equation}

\cvm
\nd
with obvious generalizations to other quantities.

\cvb
Let $M$ be a \tit{space-like} hypersurface, i.e. the induced metric is Riemannian,
with a differentiable normal $\n$ that is time-like.

In local coordinates, $(x^\al)$ and $(\x^i)$, the geometric quantities of the
space-like hypersurface $M$ are connected through the following equations

\begin{equation}\lae{2.3}
x_{ij}^\al= h_{ij}\n^\al
\end{equation}

\cvm
\nd
the so-called \tit{Gau{\ss} formula}. Here, and also in the sequel, a covariant
derivative is always a \tit{full} tensor, i.e.

\begin{equation}
x_{ij}^\al=x_{,ij}^\al-\ch ijk x_k^\al+ \cha \bet\ga\al x_i^\bet x_j^\ga.
\end{equation}

\cvm
\nd
The comma indicates ordinary partial derivatives.

In this implicit definition the \tit{second fundamental form} $(h_{ij})$ is taken
with respect to $\n$.

The second equation is the \tit{Weingarten equation}

\begin{equation}
\n_i^\al=h_i^k x_k^\al,
\end{equation}

\cvm
\nd
where we remember that $\n_i^\al$ is a full tensor.

\cvm
Finally, we have the \tit{Codazzi equation}

\begin{equation}
h_{ij;k}-h_{ik;j}=\riema\al\bet\ga\de\n^\al x_i^\bet x_j^\ga x_k^\de
\end{equation}

\cvm
\nd
and the \tit{Gau{\ss} equation}

\begin{equation}
\riem ijkl=- \{h_{ik}h_{jl}-h_{il}h_{jk}\} + \riema \al\bet\ga\de x_i^\al x_j^\bet x_k^\ga
x_l^\de.
\end{equation}

\cvm
Now, let us assume that $N$ is a globally hyperbolic Lorentzian manifold with a
\tit{compact} Cauchy surface. $N$ is then a topological product $\R[]\times \mc
S_0$, where $\mc S_0$ is a compact Riemannian manifold, and there exists a
Gaussian coordinate system
$(x^\al)$,  such that $x^0$ represents the time, the
$(x^i)_{1\le i\le n}$ are local coordinates for $\mc S_0$, where we may assume
that $\mc S_0$ is equal to the level hypersurface $\{x^0=0\}$---we don't
distinguish between $\mc S_0$ and $\{0\}\times \mc S_0$---, and such that the
Lorentzian metric takes the form

\begin{equation}\lae{2.7}
d\bar s_N^2=e^{2\psi}\{-{dx^0}^2+\s_{ij}(x^0,x)dx^idx^j\},
\end{equation}

\cvm
\nd
where $\s_{ij}$ is a Riemannian metric, $\psi$ a function on $N$, and $x$ an
abbreviation for the space-like components $(x^i)$, see \ci{GR},
\ci[p.~212]{HE}, \ci[p.~252]{GRH}, and \ci[Section~6]{cg1}.
 We also assume that
the coordinate system is \tit{future oriented}, i.e. the time coordinate $x^0$
increases on future directed curves. Hence, the \tit{contravariant} time-like
vector$(\x^\al)=(1,0,\dotsc,0)$ is future directed as is its \tit{covariant} version
$(\x_\al)=e^{2\psi}(-1,0,\dotsc,0)$.

\cvm
Let $M=\graph \fv u\so$ be a space-like hypersurface

\begin{equation}
M=\set{(x^0,x)}{x^0=u(x),\,x\in\mc S_0},
\end{equation}

\cvm
\nd
then the induced metric has the form

\begin{equation}
g_{ij}=e^{2\psi}\{-u_iu_j+\s_{ij}\}
\end{equation}

\cvm
\nd
where $\s_{ij}$ is evaluated at $(u,x)$, and its inverse $(g^{ij})=(g_{ij})^{-1}$ can
be expressed as

\begin{equation}\lae{2.10}
g^{ij}=e^{-2\psi}\{\s^{ij}+\frac{u^i}{v}\frac{u^j}{v}\},
\end{equation}

\cvm
\nd
where $(\s^{ij})=(\s_{ij})^{-1}$ and

\begin{equation}\lae{2.11}
\begin{aligned}
u^i&=\s^{ij}u_j\\
v^2&=1-\s^{ij}u_iu_j\equiv 1-\abs{Du}^2.
\end{aligned}
\end{equation}

\cvm
Hence, $\graph u$ is space-like if and only if $\abs{Du}<1$.

\cvm
We also note that

\begin{equation}\lae{2.12b}
v^{-2}=1+e^{2\psi}g^{ij}u_iu_j \equiv 1+e^{2\psi} \norm{Du}^2.
\end{equation}

\cvm
The covariant form of a normal vector of a graph looks like

\begin{equation}
(\n_\al)=\pm v^{-1}e^{\psi}(1, -u_i).
\end{equation}

\cvm
\nd
and the contravariant version is

\begin{equation}
(\n^\al)=\mp v^{-1}e^{-\psi}(1, u^i).
\end{equation}

Thus, we have

\br Let $M$ be space-like graph in a future oriented coordinate system. Then, the
contravariant future directed normal vector has the form
\begin{equation}
(\n^\al)=v^{-1}e^{-\psi}(1, u^i)
\end{equation}
and the past directed
\begin{equation}\lae{2.15}
(\n^\al)=-v^{-1}e^{-\psi}(1, u^i).
\end{equation}
\er

In the Gau{\ss} formula \re{2.3} we are free to choose the future or past directed
normal, but we stipulate that we always use the past directed normal for reasons
that we have explained in \ci{cg8}.

\cvm
Look at the component $\al=0$ in \re{2.3} and obtain in view of \re{2.15}

\begin{equation}\lae{2.16}
e^{-\psi}v^{-1}h_{ij}=-u_{ij}-\cha 000\mspace{1mu}u_iu_j- \cha 0j0
\mspace{1mu}u_i-\cha 0i0\mspace{1mu}u_j- \cha ij0.
\end{equation}

\cvm
Here, the covariant derivatives are taken with respect to the induced metric of
$M$, and

\begin{equation}
- \cha ij0=e^{-\psi}\bar h_{ij},
\end{equation}

\cvm
\nd
where $(\bar h_{ij})$ is the second fundamental form of the hypersurfaces
$\{x^0=\const\}$.

An easy calculation shows

\begin{equation}
\bar h_{ij}e^{-\psi}=-\tfrac{1}{2}\dot\s_{ij} -\dot\psi\s_{ij},
\end{equation}

\cvm
\nd
where the dot indicates differentiation with respect to $x^0$.

\cvm
Next, let us analyze under which condition a space-like hypersurface $M$ can be
written as a graph over the Cauchy hypersurface $\mc S_0$.

We first need

\cvm
\bd
Let $M$ be a  closed, space-like hypersurface in $N$. Then,
$M$ is said to be \tit{achronal}, if no two points in $M$ can be connected by a
future directed time-like curve.
\ed

\cvm
In \ci{bily78} it is proved, see also \ci[Proposition 2.5]{cg8}, 

\cvm
\bpp\lap{2.5}
Let $N$ be connected and globally hyperbolic, $\mc S_0\nobreak\ 
\su\nobreak\   N$ a compact Cauchy hypersurface, and $M\su N$ a compact,
connected space-like hypersurface of class $C^m, m\ge 1$. Then, $M=\graph \fu
u{\mc S}0$ with
$u\in C^m(\mc S_0)$ iff $M$ is achronal.
\epp

\cvb
\br\lar{2.4}
The $M_i$ are barriers for the pair $(F,f)$. Let us point out that without loss of
generality we may assume

\begin{align}
\fm 2&>f(x,\n)\qq\A\,x\in M_2,\lae{2.20b}\\
\intertext{and}
\fmo \Si &< f(x,\n)\qq\A\,x\in\Si,\lae{2.20bb}
\end{align}

\cvm
\nd for let $\h\in C^\un (\bar\Om)$ be a function with support in a small
neighbourhood of $M_1\uud M_2$---the dot should indicate that the union is
disjoint--- such that

\begin{equation}
\fv\h{M_1}>0\q\tup{and}\q \fv\h{M_2}<0
\end{equation}

\cvm
\nd and define for $\de>0$

\begin{equation}
f_\de=f+\de\h.
\end{equation}

\cvm
\nd Then, if we assume $f$ to be strictly positive with a  positive
lower bound, we have for small
$\de$

\begin{equation}
f_\de\ge \frac12 f,
\end{equation}

\cvm
\nd and the $M_i$ are barriers for $(F,f_\de)$ satisfying the strict inequalities;
since we shall derive $C^{4,\al}$ estimates independent of $\de$, we shall have
proved the existence of a solution for $f$ if we can prove it for $f_\de$.
\er

\cvb
\bl\lal{2.5b}
Let $M_i$ be barriers for $(F,f)$ satisfying the strict inequalities \re{2.20b} and
\re{2.20bb}, where $F$ is supposed to be monotone and concave. Then, they are
also barriers for the elliptic regularizations
$F_\e$ for small $\e$.
\el

\bp
In view of \rl{1.1}, we know that $\C\su \C_\e$ and $H$ is positive in $\C_\e$.
Hence, $M_2$ is certainly an upper barrier for $(F_\e,f)$ because of the
monotonicity of $F$.

Let $\Si_\e$ resp. $\Si$ be the points in $M_1$ where the principal curvatures
belong to
$\C_\e$ resp. $\C$ and assume that $\Si\ne \eS$. Suppose $M_1$ were not a
lower barrier for $(F_\e,f)$ for small $\e$, then, there exist a sequence $\e\ra
0$ and a corresponding convergent sequence $x_\e\in \Si_\e$, $x_\e\ra
x_0\in\Si$, such that

\begin{equation}
F_\e\ge f(x_\e,\n),
\end{equation}

\cvm
\nd and hence,

\begin{equation}
F\ge f(x_0,\n)
\end{equation}

\cvm
\nd contradicting \re{2.20bb}.
\ep

\cvb
\br
The condition \re{0.3} is reasonable as is evident from the Einstein equation

\begin{equation}
\bar R_{\al\bet}-\tfrac12 \bar R \msp \bar g_{\al\bet}=T_{\al\bet},
\end{equation} 

\cvm
\nd where the energy-momentum tensor $T_{\al\bet}$ is supposed to be positive
semi-definite for time-like vectors ({\it weak energy condition}, cf. \ci[p.
89]{HE}), and the relation

\begin{equation}
R=-[H^2-h_{ij}h^{ij}]+\bar R+2 \bar R_{\al\bet}\n^\al\n^\bet
\end{equation}

\cvm
\nd for the scalar curvature of a space-like hypersurface; but it would be
convenient for the approximations we have in mind, if the estimate in \re{0.3}
would be valid for all time-like vectors. 

\cvm\nd
In fact, we may assume this without loss
of generality: Let $\vt$ be a smooth real function such that

\begin{equation}
\frac{c_1}2\le \vt\q\tup{and}\q \vt(t)=t\q \A\,t\ge c_1,
\end{equation}

\cvm
\nd then, we can replace $f$ by $\vt\circ f$ and the new function satisfies our
requirements for all time-like vectors. 
 
\cvm
We therefore assume in the following that the relation \re{0.3} holds for all
time-like vectors $\n\in T_x(N)$ and all $x\in \bar\Om$.
\er

\cvb
Sometimes, we need a Riemannian reference metric, e.g. if we want to estimate
tensors. Since the Lorentzian metric can be expressed as

\begin{equation}
\bar g_{\al\bet}dx^\al dx^\bet=e^{2\psi}\{-{dx^0}^2+\s_{ij}dx^i dx^j\},
\end{equation}

\cvm
\nd
we define a Riemannian reference metric $(\tilde g_{\al\bet})$ by

\begin{equation}
\tilde g_{\al\bet}dx^\al dx^\bet=e^{2\psi}\{{dx^0}^2+\s_{ij}dx^i dx^j\}
\end{equation}

\cvm
\nd
and we abbreviate the corresponding norm of a vectorfield $\h$ by

\begin{equation}
\nnorm \h=(\tilde g_{\al\bet}\h^\al\h^\bet)^{1/2},
\end{equation}

\cvm
\nd
with similar notations for higher order tensors.

\cvb
For a space-like hypersurface $M=\graph u$ the induced metrics with respect to
$(\bar g_{\al\bet})$ resp. $( \tilde g_{\al\bet})$ are related as follows

\begin{equation}
\begin{aligned}
\tilde g_{ij}&= \tilde g_{\al\bet}x^\al_i x^\bet_j=e^{ 2\psi}[u_i u_j +\s_{ij}]\\[\cma]
&= g_{ij}+2 e^{2\psi} u_i u_j.
\end{aligned}
\end{equation}

\cvm
Thus, if $(\x^i)\in T_p(M)$ is a unit vector for $(g_{ij})$, then

\begin{equation}
\tilde g_{ij}\x^i \x^j= 1+2 e^{2\psi}\abs {u_i\x^i}^2,
\end{equation}

\cvm
\nd and we conclude for future reference

\bl\lal{2.6}
Let $M=\graph u$ be a space-like hypersurface in $N$, $p\in M$, and $\x\in
T_p(M)$ a unit vector, then

\begin{equation}
\nnorm{x^\bet_i\x^i}\le c (1+\abs{u_i\x^i})\le c \tilde v,
\end{equation}

\cvm
\nd where $\tilde v=v^{-1}$.
\el

\cvb
\section{An auxiliary curvature problem}\las{3}

\cvb
Solving the problem \re{0.2} involves two steps: first, proving a priori estimates,
and secondly, applying a method to show the existence of a solution. In a general
Lorentzian manifold the evolution method is the method of choice, but
unfortunately, one cannot prove the necessary a priori estimates during the
evolution when $F$ is the scalar curvature operator. Both the $C^1$ and
$C^2$- estimates fail for general $f=f(x,\n)$.

\cvb
Therefore, we  use the \tit{elliptic regularization} and consider the
existence problem for the operators

\begin{equation}\lae{3.1}
F_\e(\ka_i)=F(\ka_i+\e H),\qq \e>0,
\end{equation}

\cvm
\nd i.e. we solve

\begin{equation}
{F_\e}_{|_M}=f(x,\n).
\end{equation}

\cvm
Then, we prove uniform $C^{2,\al}$- estimates for the approximating solutions
$M_\e$, and finally, let $\e$ tend to zero.

\cvm
The $F_\e$---or some positive power of it---belong to a class of curvature
functions
$F$ that satisfy the following condition $(H)$: $F\in C^{2,\al}(\C)\ii
C^0(\bar\C)$, where $\C\su \R[n]$ is an open cone containing $\C_+$, $F$ is
symmetric,  monotone, i.e. $F_i>0$, homogeneous of degree 1, concave, vanishes
on
$\pa\C$,  and there exists
$\e_0=\e_0(F)>0$ such that

\begin{equation}\lae{3.3}
F_i\ge \e_0\sum_kF_k\qq \A\,1\le i\le n.
\end{equation}

\cvm
Furthermore, the set

\begin{equation}\lae{3.4}
\Lam_{\de,\ka}=\set{(\ka_i)\in\C}{0<\de\le F(\ka_i),\,\ka_i\le \ka\;\A\,1\le i\le n}
\end{equation}

\cvm
\nd
is compact.

\cvb
\br\lar{3.1}
If the original curvature function $F\in C^{2,\al}
(\C)\ii C^0(\bar\C)$ is
concave, homogeneous of degree 1, and vanishes on $\pa \C$, then, the $F_\e$
are of class $(H)$ in the cone $\C_\e$, and  satisfy \re{3.3} with $\e_0=\e$.
The set 

\begin{equation}\lae{3.4b}
\tilde \Lam_{\de,\ka}=\set{( \ka_i)\in\C_\e}{0<\de\le F_\e(\ka_i),\,\ka_i\le
\ka\;\A\,1\le i\le n}
\end{equation}

\cvm
\nd 
is   compact  for fixed $\e$. 

 If the
parameters $\ka$ and $\de$ are independent of $\e$, then the $\tilde \Lam_{\de,\ka}$
are contained in a compact subset of $\C$  uniformly in $\e$, for small
$\e$,
$0\le\e\le
\e_1(\de,\ka, F)$. 
\er

\bp
In view of the results in \rl{1.4} we only have to prove the compactness of
$\tilde \Lam_{\de,\ka}$. We shall also only consider the case when the estimates hold
uniformly in $\e$.

\cvm
Due to the concavity and homogeneity of $F_\e$ we conclude from \re{1.25}
that

\begin{equation}
F_\e(\ka_i)\le \frac1nF(1,\dots,1)(1+n\e)H.
\end{equation}

\cvm
For $(  \ka_i)\in \tilde \Lam_{\de,\ka}$ we, therefore, infer

\begin{align}
\de\le F_\e(\ka_i)\le
\frac{1+n\e}nF(1,\dots,1)H&\le(1+n\e)F(1,\dots,1)\ka,\lae{3.6}\\[\cma]
\intertext{and thus,}
\lim_{\e\ra 0}\e H&=0,\lae{3.7b}
\end{align}

\cvm
\nd uniformly in $ \tilde \Lam_{\de,\ka}$.

\cvm
Suppose $ \tilde \Lam_{\de,\ka}$ would
not stay in a  compact subset of $\C$ for small
$\e, 0<\e
\le\e_1(\de,\ka,F)$.  Then, there would exist a sequence $\e\ra 0$ and a
corresponding sequence $(  \ka^\e_i)\in \tilde \Lam_{\de,\ka}$ converging to
a point
$(\ka_i)\in
\pa\C$, which is impossible in view of \re{3.6},
\re{3.7b}, and the continuity of $F$ in $\bar \C$.
\ep

\cvb
To prove the existence of hypersurfaces of prescribed curvature $F$ for $F\in
(H)$ we look at the evolution problem

\begin{equation}
\begin{aligned}\lae{3.8}
\dot x&=(F-f)\n,\\[\cma]
x(0)&=x_0,
\end{aligned}
\end{equation}

\cvm
\nd where $\n$ is the past-directed normal of the flow hypersurfaces $M(t)$,
$F$ the curvature evaluated at $M(t)$, $x=x(t)$ an embedding and $x_0$ an
embedding of an initial hypersurface $M_0$, which we choose to be the upper
barrier $M_2$.

Since $F$ is an elliptic operator,  short-time existence, and hence, existence in a
maximal time interval $[0,T^*)$ is guaranteed. If we are able to prove uniform a
priori estimates in $C^{2,\al}$, long-time existence and convergence to a
stationary solution will follow immediately.

\cvb
But before we prove the a priori estimates, we want to show how the metric, the
second fundamental form, and the normal vector of the hypersurfaces $M(t)$
evolve. All time derivatives are
\tit{total} derivatives. The proofs are identical to those of the corresponding
results in a Riemannian setting, cf. \ci[Section 3]{cg2} and \ci[Section 4]{cg8},
and will be omitted.

\cvb
\bl[Evolution of the metric]
The metric $g_{ij}$ of $M(t)$ satisfies the evolution equation
\begin{equation}
\dot g_{ij}=2( F- f)h_{ij}.
\end{equation}
\el

\cvb
\bl[Evolution of the normal]
The normal vector evolves according to
\begin{equation}\lae{3.7}
\dot \n=\nabla_M( F- f)=g^{ij}( F- f)_i x_j.
\end{equation}
\el

\cvb
\bl[Evolution of the second fundamental form]
The second fundamental form evolves according to

\begin{equation}
\dot h_i^j=( F- f)_i^j- ( F- f) h_i^k h_k^j - ( F- f) \riema
\al\bet\ga\de\n^\al x_i^\bet \n^\ga x_k^\de g^{kj}
\end{equation}
and
\begin{equation}
\dot h_{ij}=( F- f)_{ij}+ ( F- f) h_i^k h_{kj}- ( F- f) \riema
\al\bet\ga\de\n^\al x_i^\bet \n^\ga x_j^\de.
\end{equation}
\el

\cvb
\bl[Evolution of $( F- f)$]
The term $( F- f)$ evolves according to the equation

\begin{multline}\lae{3.10}
{( F- f)}^\prime- F^{ij}( F- f)_{ij}=\msp[3]- 
F^{ij}h_{ik}h_j^k ( F- f)
 - f_\al\n^\al ( F- f)\\-f_{\n^\al}x^\al_i(F-f)_jg^{ij}
- F^{ij}\riema \al\bet\ga\de\n^\al x_i^\bet \n^\ga x_j^\de ( F- f),
\end{multline}
\el

\cvb
From \re{3.8} we deduce with the help of the Ricci identities a parabolic equation
for the second fundamental form

\bl\lal{3.6}
The mixed tensor $h_i^j$ satisfies the parabolic equation

\begin{equation}\lae{3.14}
\begin{aligned}
\dot h_i^j- &F^{kl}h_{i;kl}^j\\
&=-  F^{kl}h_{rk}h_l^rh_i^j
+f h_i^kh_k^j\\
&\hp{+}- f_{\al\bet} x_i^\al x_k^\bet g^{kj}-  f_\al\n^\al
h_i^j-f_{\al\n^\bet}(x^\al_i x^\bet_kh^{kj}+x^\al_l
x^\bet_k h^{k}_i\, g^{lj})\\ &\hp{=}
-f_{\n^\al\n^\bet}x^\al_lx^\bet_kh^k_ih^{lj}-f_{\n^\bet} x^\bet_k h^k_{i;l}\,g^{lj} 
-f_{\n^\al}\n^\al h^k_i h^j_k\\ &\hp{=}+
F^{kl,rs}h_{kl;i}h_{rs;}^{\hphantom{rs;}j}+2 F^{kl}\riema
\al\bet\ga\de x_m^\al x_i ^\bet x_k^\ga x_r^\de h_l^m g^{rj}\\
&\hp{=}- F^{kl}\riema \al\bet\ga\de x_m^\al x_k ^\bet x_r^\ga x_l^\de
h_i^m g^{rj}- F^{kl}\riema \al\bet\ga\de x_m^\al x_k ^\bet x_i^\ga x_l^\de h^{mj} \\
&\hp{=}- F^{kl}\riema \al\bet\ga\de\n^\al x_k^\bet\n^\ga x_l^\de h_i^j+ f
\riema \al\bet\ga\de\n^\al x_i^\bet\n^\ga x_m^\de g^{mj}\\
&\hp{=}+ F^{kl}\bar R_{\al\bet\ga\de;\e}\{\n^\al x_k^\bet x_l^\ga x_i^\de
x_m^\e g^{mj}+\n^\al x_i^\bet x_k^\ga x_m^\de x_l^\e g^{mj}\}.
\end{aligned}
\end{equation}
\el

\cvb
The proof is identical to that of the corresponding result in the Riemannian case,
cf. \ci[Lemma 7.1 and Lemma 7.2]{cg2}; the only difference is that $f$ now also
depends on $\n$.

\cvb
\br\lar{3.7}
In view of the maximum principle, we immediately deduce from \re{3.10} that the
term $( F- f)$ has a sign during the evolution if it has one at the beginning,
i.e., if the starting hypersurface $M_0$ is the upper barrier $M_2$, then
$( F- f)$ is non-negative

\begin{equation}\lae{3.15}
F\ge f.
\end{equation}
\er

\cvb
\section{Lower order estimates for the auxiliary solutions}\las{4}

\cvb
Since the two boundary components $M_1, M_2$ of $\pa\Om$ are  space-like,
achronal hypersurfaces, they can be written as graphs over the Cauchy
hypersurface $\so$, $M_i=\graph u_i$, $i=1,2$, and we have
\begin{equation}\lae{4.1}
u_1\le u_2,
\end{equation}
for $M_1$ should lie in the past of $M_2$, and the enclosed domain is supposed to
be connected. 

Let us look at the evolution equation \re{3.8} with initial hypersurface $M_0$
equal to $M_2$ defined on  a maximal time interval
$I=[0,T^*),T^*\le
\un$. 
Since the initial hypersurface is a graph over $\so$, we can write

\begin{equation}
M(t)=\graph\fu{u(t)}S0\q \A\,t\in I,
\end{equation}
where $u$ is defined in the cylinder $Q_{T^*}=I\times \so$. We then deduce from
\re{3.8}, looking at the component $\al=0$, that $u$ satisfies a parabolic
equation of the form
\begin{equation}\lae{4.3}
\dot u=-e^{-\psi}v^{-1}(F- f),
\end{equation}
where we  use the notations in \rs{2}, and where we emphasize that the time
derivative is a total derivative, i.e.

\begin{equation}\lae{4.4}
\dot u=\pde ut+u_i\dot x^i.
\end{equation}

Since the past directed normal can be expressed as

\begin{equation}
(\n^\al)=-e^{-\psi}v^{-1}(1,u^i),
\end{equation}

\nd
we conclude from \re{3.8}, \re{4.3}, and \re{4.4}

\begin{equation}\lae{4.6}
\pde ut=-e^{-\psi}v(F- f).
\end{equation}

\cvm
Thus, $\pde ut$ is non-positive in view of \rr{3.7}.

\cvb
Next, let us state our first a priori estimate

\cvm
\bl\lal{4.1}
Suppose that the boundary components act as barriers for $(F,f)$, then  the flow hypersurfaces stay in
$\bar
\Om$ during
the evolution.
\el

\cvb
The proof is identical to that of the corresponding result in \ci[Lemma 4.1]{cg8}.

\cvb
For the $C^1$- estimate the term $\tilde v=v^{-1}$ is of great importance. It
satisfies the following evolution equation

\cvb
\bl[Evolution of $\tilde v$]\lal{4.2}
Consider the flow \re{3.8} in the distinguished coordinate system associated
with $\so$. Then, $\tilde v$ satisfies the evolution equation

\begin{equation}\lae{4.7}
\begin{aligned}
\dot{\tilde v}- F^{ij}\tilde v_{ij}=&- F^{ij}h_{ik}h_j^k\tilde v
-f\h_{\al\bet}\n^\al\n^\bet\\
&-2 F^{ij}h_j^k x_i^\al x_k^\bet \h_{\al\bet}- F^{ij}\h_{\al\bet\ga}x_i^\bet
x_j^\ga\n^\al\\
&- F^{ij}\riema \al\bet\ga\de\n^\al x_i^\bet x_k^\ga x_j^\de\h_\e x_l^\e g^{kl}\\
&- f_\bet x_i^\bet x_k^\al \h_\al g^{ik} -f_{\n^\bet}x^\bet_k h^{ik}x^\al_i\h_\al,
\end{aligned}
\end{equation}
where $\h$ is the covariant vector field $(\h_\al)=e^{\psi}(-1,0,\dotsc,0)$.
\el

\cvb
The proof uses the relation

\begin{equation}\lae{4.8}
\tilde v=\h_\al \n^\al
\end{equation}

\cvm
\nd and
is identical to that of \ci[Lemma 4.4]{cg8} having in
mind that  presently $f$  also depends on $\n$.

\cvb
\bl\lal{4.3}
Let $M(t)=\graph u(t)$ be the flow hypersurfaces, then, we have

\begin{equation}\lae{4.9}
\begin{aligned}
\dot u-F^{ij}u_{ij}=e^{-\psi}\tilde v f&+\cha 000\, F^{ij}u_i u_j \\
&+2\msp
F^{ij}\cha 0i0\,u_j +F^{ij}\cha ij0,
\end{aligned}
\end{equation}

\cvm
\nd where all covariant derivatives a taken with respect to the induced metric of
the flow hypersurfaces, and the time derivative $\dot u$ is the total time
derivative, i.e. it is given by \re{4.4}.
\el

\cvb
\bp
We use the relation \re{4.3} together with \re{2.16}.
\ep

\cvb
As an immediate consequence we obtain

\cvb
\bl\lal{4.4}
The composite function

\begin{equation}
\f=e^{\m e^{\lam u}}
\end{equation}

\cvm
\nd where $\m,\lam$ are constants, satisfies the equation

\begin{equation}\lae{4.10}
\begin{aligned}
\dot\f -F^{ij}\f_{ij}=&fe^{-\psi} \tilde v\msp \m \lam \msp[2] e^{\lam u} \,\f + F^{ij}
u_i u_j \,\cha000\,\m \lam\,e^{\lam u}\f\\[\cma]
&+2\msp F^{ij} u_i \cha 0j0\,\m \lam\msp[2] e^{\lam u}\,\f+ F^{ij}\cha
ij0\,\m\lam\msp[2] e^{\lam u}\,\f\\[\cma]
&-[1+\m\msp[2] e^{\lam u}] F^{ij} u_i u_j\,\m \lam^2\msp[2] e^{\lam u}\,\f.
\end{aligned}
\end{equation}
\el

\cvb
Before we can prove the $C^1$- estimates we need two more lemmata.

\cvb
\bl\lal{4.5}
There is a constant $c=c(\Om)$ such that for any positive function $0<\e=\e(x)$
on $\so$ and any hypersurface $M(t)$ of the flow we have

\begin{align}
\nnorm{\n}&\le c \msp \tilde v,\\[\cma]
g^{ij}&\le c\msp \tilde v^2\msp \s^{ij},\lae{4.12}\\[\cma]
F^{ij}&\le F^{kl} g_{kl}\msp g^{ij},\lae{4.13}
\end{align}

\begin{equation}\lae{4.14}
\begin{aligned}
\abs{F^{ij}h^k_j x^\al_i x^\bet_k \msp\h_{\al\bet}}\le \frac\e 2F^{ij}h^k_i
h_{kj}\msp
\tilde v + \frac c{2\e} F^{ij}g_{ij}\msp  \tilde v^3,
\end{aligned}
\end{equation}

\begin{equation}\lae{4.15}
\abs{F^{ij}\h_{\al\bet\ga} x^\bet_i x^\ga_j \n^\al}\le c\msp \tilde v^3 F^{ij}g_{ij},
\end{equation}

\nd and
\begin{equation}\lae{4.16}
\abs{F^{ij}\riema \al\bet\ga\de \n^\al x^\bet_i x^\ga_k x^\de_j\h_\e x^\e_l g^{kl}}\le
c\msp
\tilde v^3 F^{ij}g_{ij}.
\end{equation}
\el

\cvb
\bp
(i) The first three inequalities are obvious.

\cvm
(ii) \re{4.14} follows from the generalized Schwarz inequality combined with
\re{4.12} and \re{4.13}.

\cvm
(iii) \re{4.15} is a direct consequence of \re{4.12} and \re{4.13}.

\cvm
(iv) The proof of \re{4.16} is a bit more complicated and uses the symmetry
properties of the Riemann curvature tensor.

Let

\begin{equation}\lae{4.17}
a_{ij}=\riema \al\bet\ga\de \n^\al x^\bet_i x^\ga_k x^\de_j \h_\e x^\e_l g^{kl}.
\end{equation}

\cvm
We shall show that the symmetrization of $a_{ij}$ satisfies

\begin{equation}\lae{4.18}
-c\msp \tilde v^3 g_{ij}\le \frac12 (a_{ij}+a_{ji})\le c\msp \tilde v^3 g_{ij}
\end{equation}

\cvm
\nd with a uniform constant $c=c(\Om)$, which in turn yields \re{4.16}.

\cvm
Let $p\in M(t)$ be  arbitrary,  $(x^\al)$ be the special Gaussian coordinate of
$N$, and $(\x^i)$ local coordinates around $p$ such that

\begin{equation}
x^\al_i=
\begin{cases}
 u_i, &\al= 0,\\
\de^k_i , &\al=k.
\end{cases}
\end{equation}

\cvm
We also note that all indices are raised with respect to $g^{ij}$ with the exception
of the contravariant vector

\begin{equation}
\check u^i =\s^{ij}u_j.
\end{equation}

\cvm
We point out that

\begin{align}
\norm{Du}^2&=g^{ij}u_i u_j=e^{-2\psi} \tilde v^2 \s^{ij} u_i u_j,\\[\cma]
\tilde v^2&=1+e^{2\psi} \norm{Du}^2,\lae{4.22}\\[\cma]
(\n^\al)&=- \tilde v(1, \check u^i) e^{-\psi},\\
\intertext{and}
\h_\e x^\e_l g^{kl}&=-e^\psi u^k.\lae{4.24}
\end{align}

\cvm
First, let us observe that in view of \re{4.24} and the symmetry properties of
the Riemann curvature tensor we have

\begin{equation}\lae{4.25}
a_{ij}u^j=0.
\end{equation}

\cvm
Next, we shall expand the right-hand side of \re{4.17} explicitly.

\begin{equation}\lae{4.26}
\begin{aligned}
a_{ij}&=  \riema 0i0j \tilde  v\norm{Du}^2 + \riema 0ik0 \tilde v u_j u^k
+\riema 0ikj \tilde  vu^k\\[\cma]
&\hp{=}+ \riema l0k0 \tilde  vu^k \check u^l u_i u_j + \riema l00j \tilde v \check
u^l u_i \norm{Du}^2\\[\cma]
&\hp{=} +\riema l0kj \tilde v u^k\check u^l u_i +\riema li0j \tilde v \check
u^l\norm{Du}^2\\[\cma]
&\hp{=} +\riema lik0 \tilde v u^k \check u^lu_j +\riema likj \tilde v u^k \check
u^l
\end{aligned}
\end{equation}

\cvb
To prove the estimate \re{4.18}, we may assume that $Du\ne 0$. Let $e_i$,
$1\le i\le n$ be an orthonormal base of $T_p(M(t))$ such that

\begin{equation}
e_1=\frac{Du}{\norm{Du}},
\end{equation}

\cvm
\nd then, for $2\le k\le n$, the $e_k$ are also orthonormal with respect to the
metric $e^{2\psi} \s_{ij}$, and it is also valid that

\begin{equation}
\s_{ij} \check u^i e^j_k=0\qq\A\;2\le k\le n,
\end{equation}

\cvm
\nd where $e_k=(e^i_k)$.

\cvm
For $2\le r,s\le n$ we deduce from \re{4.26}

\begin{equation}
\begin{aligned}
a_{ij}e^i_r e^j_s &= \riema 0i0j \tilde v\norm{Du}^2 e^i_r e^j_s +\riema 0ikj
\tilde v u^k e^i_r e^j_s\\[\cma]
&\hp{=} +\riema li0j\tilde v \check u^l\norm{Du}^2 e^i_r e^j_s+ \riema likj \tilde
v u^k \check u^l e^i_r e^j_s
\end{aligned}
\end{equation}

\cvm
\nd and hence,

\begin{equation}
\abs{a_{ij}e^i_r e^j_s}\le c \msp\tilde v^3\qq\A\; 2\le r,s\le n.
\end{equation}

\cvm
It remains to estimate $a_{ij}e^i_1 e^j_r$ for $2\le r\le n$, because of
\re{4.25}.

\cvm
We deduce from \re{4.26}

\begin{equation}
a_{ij} e^i_1 e^j_r=\riema 0i0j \tilde v\norm{Du}^2 \tilde v^{-2} e^i_1 e^j_r+
\riema 0ikj \tilde v^{-1} u^k e^i_1 e^j_r,
\end{equation}

\cvm
\nd where we used the symmetry properties of the Riemann curvature tensor.

Hence, we conclude

\begin{equation}
\abs{a_{ij} e^i_1 e^j_r}\le c\msp \tilde v^2\qq\A\; 2\le r\le n,
\end{equation}

\cvm
\nd and the relation \re{4.18} is proved.
\ep

\cvb
\bl\lal{4.6}
Let $M\su \bar \Om$ be a graph over $\so$, $M=\graph u$, and
$\e=\e(x)$ a function given in $\so$, $0<\e <\frac12$. Let $\f$ be
defined through

\begin{equation}\lae{4.33}
\f=e^{\m e^{\lam u}},
\end{equation}

\cvm
\nd where $0<\m$ and $\lam<0$. Then, there exists $c=c(\Om)$ such that

\begin{equation}
\begin{aligned}
2\abs{F^{ij} \tilde v_i \f_j}&\le c\msp F^{ij}g_{ij} \tilde v^3 \abs{\lam} \m e^{\lam u}
\f +(1-2\e) F^{ij} h^k_i h_{kj} \tilde v \f\\[\cma]
&\hp{\le} +\frac1{1-2\e} F^{ij} u_i u_j \m^2 \lam^2 e^{2\lam u} \tilde v \f.
\end{aligned}
\end{equation}
\el

\cvb
\bp
Since $ \tilde v=\h_\al \n^\al$, we have

\begin{equation}\lae{4.35}
\begin{aligned}
\tilde v_i&= \h_{\al\bet}\n^\al x^\bet_i + \h_\al h^k_i x^\al_k\\[\cma]
&= \h_{\al\bet} \n^\al x^\bet_i -e^\psi h^k_i u_k.
\end{aligned}
\end{equation}

\cvm
Thus, we derive

\begin{equation}
\begin{aligned}
2 \abs{F^{ij} \tilde v_i\f_j}&= 2\abs{F^{ij} \tilde v_i u_j} \abs\lam \m e^{\lam u}
\f\\[\cma]
&\le c\msp  F^{ij} g_{ij} \tilde v^3 \abs\lam \m e^{\lam u} \f +2 e^\psi \abs{F^{ij}
h^k_i u_k u_j}\abs\lam \m e^{\lam u} \f.
\end{aligned}
\end{equation}

\cvm
The last term of the preceding inequality can be estimated by 

\begin{equation}
(1-2\e) F^{ij} h^k_i h_{kj} \tilde v\f +\frac1{1-2\e} \tilde v^{-1} e^{2\psi}
\norm{Du}^2 F^{ij} u_i u_j \m^2 \lam^2 e^{2\lam u} \f
\end{equation}

\cvm
\nd and we obtain the desired estimate in view of \re{4.22}.
\ep

\cvb
Applying \rl{4.5} to the evolution equation for $ \tilde v$ we conclude

\cvb
\bl\lal{4.7}
There exists a constant $c=c(\Om)$ such that for any function $\e$,
$0<\e=\e(x)<1$, defined on $\so$ the term $ \tilde v$ satisfies an evolution
inequality of the form

\begin{equation}\lae{4.39b}
\begin{aligned}
\dot {\tilde v} -F^{ij} \tilde v_{ij}&\le -(1-\e) F^{ij} h^k_i h_{kj} \tilde v -f
\h_{\al\bet} \n^\al \n^\bet\\[\cma]
&\hp{\le} +\frac{c}\e F^{ij} g_{ij} \tilde v^3 +c\msp \nnorm{f_\bet} \tilde v^2+
f_{\n^\bet} x^\bet_l h^{kl} u_k e^\psi.
\end{aligned}
\end{equation}
\el

\cvb
We are now ready to prove the uniform boundedness of $ \tilde v$.

\cvb
\bpp\lap{4.8}
Assume that there are positive constants $c_i$, $1\le i\le 3$, such that for any
$x\in \Om$ and any past directed time-like vector $\n$ there holds

\begin{align}
-c_1&\le f(x,\n),\lae{4.39}\\[\cma]
\nnorm{f_\bet(x,\n)}&\le c_2 (1+\nnorm{\n}),\lae{4.40}\\
\intertext{and}
\nnorm{f_{\n^\bet}(x,\n)}&\le c_3.\lae{4.41}
\end{align}

\cvm
Then, the term $ \tilde v$ remains uniformly bounded during the evolution

\begin{equation}
\tilde v\le c=c(\Om, c_1, c_2, c_3, \e_0),
\end{equation}

\cvm
\nd where $\e_0$ is the constant in \re{3.3}. Here, and in the following, the
reference that a constant depends on $\Om$ also means that it depends on the
barriers and geometric quantities of the ambient space restricted to $\Om$.
\epp

\cvb
\bp
We proceed similar as in \ci[Proposition 3.7]{cg7} and show that the function

\begin{equation}\lae{4.44}
w= \tilde v\f,
\end{equation}

\cvm
\nd $\f$ as in \re{4.33}, is uniformly bounded, if we choose

\begin{equation}\lae{4.45}
0<\m<1\q \tup{and}\q \lam<<-1,
\end{equation}

\cvm
\nd appropriately,  and assume furthermore, without loss of generality, that $u\le
-1$, for otherwise replace $u$ by $(u-c)$, $c$ large, in the definition of $\f$.

\cvb
With the help of  the lemmata \ref{L:4.4}, \ref{L:4.6}, and \ref{L:4.7} we derive
from the relation

\begin{equation}
\begin{aligned}
\dot w - F^{ij} w_{ij}=[\dot{ \tilde v}- F^{ij} \tilde v_{ij}] \f+
[\dot\f-F^{ij}\f_{ij}] \tilde v-2 F^{ij} \tilde v_i \f_j
\end{aligned}
\end{equation}

\cvm
\nd the parabolic inequality

\begin{equation}\lae{4.47b}
\begin{aligned}
\dot w -F^{ij} w_{ij}&\le -\e\msp[2] F^{ij} h^k_i h_{kj} \tilde v \f +
c[\e^{-1}+\abs\lam \m e^{\lam u}]F^{ij} g_{ij} \tilde v^3 \f\\[\cma]
&\hp{\le} +[\frac1{1-2\e}-1] F^{ij} u_i u_j \m^2\lam^2 e^{2\lam u} \tilde v
\f \\[\cma]
&\hp{\le}-F^{ij} u_i u_j \m\lam^2 e^{\lam u} \tilde v \f \\[\cma]
&\hp{\le} +f
[-\h_{\al\bet}\n^\al \n^\bet+e^{-\psi}
\m
\lam e^{\lam u} \tilde v^2] \f\\[\cma]
&\hp{\le} + c\msp[2] \nnorm{f_\bet} \tilde v^2 \f +f_{\n^\bet} x^\bet_l h^{kl} u_k
e^\psi
\f,
\end{aligned}
\end{equation}

\cvm
\nd where we have chosen the same function $\e=\e(x)$ in \rl{4.6} resp.
\rl{4.7}.

\cvb
Setting $\e=e^{-\lam u}$ and using \rl{2.6}, the assumption \re{3.3}, which can be
rewritten as

\begin{equation}\lae{4.47}
F^{ij}\ge \e_0 F^{kl}g_{kl} g^{ij},
\end{equation}

\cvm
\nd as well as the assumptions \re{4.39}, \re{4.40}, and \re{4.41}, and
observing, furthermore, that in view of the concavity and homogeneity  of $F$

\begin{equation}\lae{4.48}
F^{ij}g_{ij}\ge F(1,\dotsc,1) >0,
\end{equation}

\cvm
\nd we conclude

\begin{equation}\lae{4.50}
\begin{aligned}
\dot w -F^{ij} w_{ij}&\le -\frac12 F^{ij}h^k_i h_{kj}  e^{-\lam u}\tilde v\f + c
\abs\lam \m e^{\lam u} F^{ij}g_{ij} \tilde v^3 \f\\[\cma]
&\hp{\le } +\frac2{1-2\e} F^{ij} u_i u_j \m^2\lam^2 e^{\lam u} \tilde v
\f -F^{ij} u_i u_j \m\lam^2 e^{\lam u} \tilde v \f \\[\cma]
&\hp{\le}+ c \msp c_1 \m \abs\lam e^{\lam u} \tilde v^2 \f + c\msp c_2 \tilde v^3 \f
+c\msp c^2_3 \e_0^{-1} e^{\lam u}\tilde v^3 \f,
\end{aligned}
\end{equation}

\cvm
\nd  where $\abs\lam$ is chosen so large that

\begin{equation}
e^{-\lam u}\le \frac14.
\end{equation}

\cvm
Choosing, furthermore,

\begin{equation}
\m=\frac18,
\end{equation}

\cvm
\nd we see that the terms involving $F^{ij} u_i u_j$ add up to a dominating
negative quantity that can be estimated from above by

\begin{equation}
-\frac1{16} F^{ij} u_i u_j \lam^2 e^{\lam u} \tilde v \f\le -\frac{\e_0}{16} F^{kl}
g_{kl} \norm{Du}^2 \lam^2 e^{\lam u} \tilde v \f,
\end{equation}

\cvm
\nd in view of \re{4.47}.

\cvm
$\norm{Du}^2$ is of the order $ \tilde v^2$ for large $ \tilde v$, hence, the
parabolic maximum principle yields a uniform estimate for $w$ if $\abs\lam$ is
chosen large enough.
\ep

\cvb
\section{$C^2$- estimates for the auxiliary solutions}\las{5}

\cvb
We want to prove that the principal curvatures of the flow hypersurfaces are
uniformly bounded.

\cvb
\bpp\lap{5.1}
Let $M(t)$, $0\le t<T^*$, be solutions of the evolution problem \re{3.8} with
$M(0)=M_2$, $F\in (H)$, and $f\in C^{2,\al}$ strictly positive,

\begin{equation}\lae{5.1}
0<c_0\le f.
\end{equation}

\cvm
\nd Then, the principal curvatures of the flow hypersurfaces are uniformly
bounded provided the $M(t)$ are uniformly space-like, i.e. uniform $C^1$-
estimates are valid.
\epp

\cvb
\bp
As we have already mentioned in \rr{3.7}, we know that 

\begin{equation}\lae{5.2}
0<c_0\le f\le F
\end{equation}

\cvm
\nd during the evolution, thus, it is sufficient  to estimate the principal
curvatures from above.

Let $\f$  be defined  by

\begin{equation}\lae{5.3}
\f=\sup\set{{h_{ij}\h^i\h^j}}{{\norm\h=1}}.
\end{equation}
\cvm
\nd
We claim that
$\f$ is uniformly bounded.

Let $0<T<T^*$, and $x_0=x_0(t_0)$, with $ 0<t_0\le T$, be a point in $M(t_0)$ such
that

\begin{equation}
\sup_{M_0}\f<\sup\set {\sup_{M(t)} \f}{0<t\le T}=\f(x_0).
\end{equation}

We then introduce a Riemannian normal coordinate system $(\x^i)$ at $x_0\in
M(t_0)$ such that at $x_0=x(t_0,\x_0)$ we have
\begin{equation}
g_{ij}=\de_{ij}\q \tup{and}\q \f=h_n^n.
\end{equation}

Let $\tilde \h=(\tilde \h^i)$ be the contravariant vector field defined by
\begin{equation}
\tilde \h=(0,\dotsc,0,1),
\end{equation}
and set
\begin{equation}
\tilde \f=\frac{h_{ij}\tilde \h^i\tilde \h^j}{g_{ij}\tilde \h^i\tilde \h^j}\raise 2pt
\hbox{.}
\end{equation}

$\tilde \f$ is well defined in neighbourhood of $(t_0,\x_0)$, and $\tilde \f$
assumes its maximum at $(t_0,\x_0)$. Moreover, at $(t_0,\x_0)$ we have
\begin{equation}
\dot{\tilde \f}=\dot h_n^n,
\end{equation}
and the spatial derivatives do also coincide; in short, at $(t_0,\x_0)$ $\tilde \f$
satisfies the same differential equation \re{3.14} as $h_n^n$. For the sake of
greater clarity, let us therefore treat $h_n^n$ like a scalar and pretend that
$\f=h_n^n$. 

\cvb
At $(t_0,\x_0)$ we have $\dot\f\ge 0$, and, in view of the maximum principle, we
deduce from \rl{3.6}

\begin{equation}\lae{5.9}
\begin{aligned}
0&\le -\e_0\msp F^{ij} g_{ij} \abs A^2 h^n_n + f\abs{h^n_n}^2 + c\msp
F^{ij}g_{ij} (h^n_n+1)\\[\cma] &\hp{\le}+c(1+\abs A^2)
(1+f+\nnorm{Df}+\nnorm{D^2f}),
\end{aligned}
\end{equation}

\cvm
\nd where we used the concavity of $F$, the Codazzi equations,  \re{4.47}, and
where

\begin{equation}
\abs A^2= g^{ij}h^k_ih_{kj}.
\end{equation}

\cvb
  Thus, $\f$ is uniformly bounded in view of \re{4.48}.
\ep

\cvb
\section{Convergence to a stationary solution}\las{6}

\cvb

We shall show that the solution of the evolution problem \re{3.8} exists for all
time, and that it converges to a stationary solution.

\cvb
\bpp
The solutions $M(t)=\graph u(t)$ of the evolution problem \re{3.8} with $F\in
(H)$, and $M(0)=M_2$ exist for all time and converge to a stationary solution
provided $f\in C^{2,\al}$ satisfies the conditions \re{4.40}, \re{4.41}, and
\re{5.1}.
\epp

\cvb
\bp
  Let us look at the scalar
version of the flow as in
\re{4.6}

\begin{equation}\lae{6.1}
\pde ut=-e^{-\psi}v(F- f).
\end{equation}

\cvm
This is  a scalar parabolic differential equation defined on the cylinder

\begin{equation}
Q_{T^*}=[0,T^*)\times \so
\end{equation}

\cvm
\nd 
with initial value $u(0)=u_2\in C^{4,\al}(\so)$. In view of the a priori estimates,
which we have established in the preceding sections, we know that

\begin{equation}
{\abs u}_\low{2,0,\so}\le c
\end{equation}
and
\begin{equation}
F\,\tup{is uniformly elliptic in}\,u
\end{equation}

\cvm
\nd
independent of $t$ due to the definition of the class $(H)$. Thus, we can apply
the known regularity results, see e.g. \ci[Chapter 5.5]{nk}, where even more
general operators are considered,  to conclude that uniform
$C^{2,\al}$-estimates are valid, leading further to uniform $C^{4,\al}$-estimates
due to the regularity results for linear operators.

\cvm
Therefore, the maximal time interval is unbounded, i.e. $T^*=\un$.

\cvm
Now, integrating \re{6.1} with respect to $t$, and observing that the right-hand
side is non-positive, yields

\begin{equation}
u(0,x)-u(t,x)=\int_0^te^{-\psi}v(F- f)\ge c\int_0^t(F- f),
\end{equation}
i.e.,
\begin{equation}
\int_0^\un \abs{F- f}<\un\qq\A\msp x\in \so.
\end{equation}

\cvm

Hence, for any $x\in\so$ there is a sequence $t_k\rightarrow \un$ such that
$(F- f)\rightarrow 0$.

On the other hand, $u(\cdot,x)$ is monotone decreasing and therefore
\begin{equation}
\lim_{t\rightarrow \un}u(t,x)=\tilde u(x)
\end{equation}
exists and is of class $C^{4,\al}(\so)$ in view of the a priori estimates. We, finally,
conclude that $\tilde u$ is a stationary solution, and that

\begin{equation}
\lim_{t\rightarrow \un}(F- f)=0.
\end{equation}
\ep

\cvb
An immediate consequence of the results we have proved so far---cf. especially
\rl{2.5b} and \rr{3.1}--- is the following theorem which is of independent interest.

\cvb
\bt
Let $F\in C^{2,\al}(\C)\ii C^0(\bar\C)$ be a concave curvature function
vanishing on $\pa\C$ and homogeneous of degree 1. Let $f=f(x,\n)$ of class
$C^{2,\al}$ satisfy the conditions
\re{4.40},
\re{4.41}, and \re{5.1}, and suppose that the boundary  components $M_i$ act
as barriers for $(F,f)$, then, there exists an admissible hypersurface $M=\graph
u$,
$u\in C^{4,\al}(\bar\so)$, solving

\begin{equation}
\fv {F_\e}M=f(x,\n)
\end{equation}

\cvm
\nd for small $\e>0$.
\et

\cvb
\section{Stationary approximations}\las{7}

\cvb
We want to solve the equation

\begin{equation}\lae{7.1}
\fv{H_2}M= f(x,\n),
\end{equation}

\cvm
\nd where $f$ satisfies the conditions of \rt{0.2}. The curvature function
$F=H^{\frac12}_2$ is concave and the elliptic regularization $F_\e$   of class
$(H)$,  cf. \re{3.1} and \rr{3.1}.

\cvm
Thus, we would like to apply the preceding existence result to find hypersurfaces
$M_\e\su\bar\Om$ such that

\begin{equation}\lae{7.2}
\fv{F_\e}{M_\e}=f^{\frac12}.
\end{equation}

\cvm
But, unfortunately, the derivatives $f_\bet$ grow quadratically in $\nnorm\n$
contrary to the assumption \re{4.40} in \rp{4.8}.

\cvm
Therefore, we define  a smooth cut-off function $\tht\in C^\un
(\R[]_+)$, $0<\tht\le 2 k$, where  $k\ge k_0>1$ is to be determined later,  by

\begin{equation}\lae{7.3}
\tht(t)=
\begin{cases}
 t, &0\le t\le k,\\
2k , &2k\le t,
\end{cases}
\end{equation}
\nd such that

\begin{equation}
0\le\dot\tht\le 4,
\end{equation}

\cvm
\nd and consider the problem

\begin{equation}\lae{7.5}
\fv{F_\e}{M_\e}= \tilde f(x, \tilde  \n),
\end{equation}

\cvm
\nd where for a space-like hypersurface $M=\graph u$ with past directed
normal vector $\n$ as in \re{2.15}, we set

\begin{align}
\tilde \n&=\tht( \tilde v) \tilde v^{-1}\n\\
\intertext{and}
\tilde f(x, \tilde \n)&=f^{\frac12}(x, \tilde \n).
\end{align}

\cvm
Then,

\begin{equation}
\nnorm{\tilde \n}\le c\msp k,
\end{equation}

\cvm
\nd so that the assumptions in \rp{4.8} are certainly satisfied.

\cvm
The constant $k_0$ should be so large that $ \tilde \n=\n$ in case of the barriers
$M_i$, $i=1,2$. 

\cvm
If we now start with the evolution equation

\begin{equation}
\dot x=(F_\e- \tilde f)\n,
\end{equation}

\cvm
\nd
then,  the $M_i$ are barriers for $(F_\e, \tilde f)$ for small $\e$, cf. \rl{2.5b}
and we conclude

\cvm
\bl\lal{7.1b}
The flow hypersurfaces $M_\e(t)=\graph u_\e$ stay in $\bar \Om$ during the
evolution if
$\e$ is small
$0<\e\le \bar\e(\Om)$.
\el

\cvb
\br
When we consider the elliptic regularizations $F_\e$, we would like to generalize
the meaning of {\it admissible} hypersurface by calling a hypersurface admissible
if the tensor $h_{ij}+\e H g_{ij}$ is admissible, i.e. if its eigenvalues belong to
$\C_2$.
\er

\cvb
Next, let us consider   the evolution equations for $ \tilde v$ and $h^j_i$
which look slightly different: In
\re{4.7} the term involving $f_{\n_\bet}$ has to be replaced by

\begin{equation}
- \tilde f_{ \tilde \n_\bet} [\tht( \tilde v) \tilde v^{-1} \n^\bet_i+\dot\tht \tilde v_i
\tilde v^{-1}\n^\bet - \tht \tilde v^{-2} \tilde v_i\n^\bet] x^\al_k g^{ik}\h_\al.
\end{equation}

\cvm
But in view of \re{4.35}, the additional terms do not cause any new problems in
the proof of \rp{4.8}, and hence, the uniform $C^1$- estimates are still valid for
the modified evolution problem, where the estimates depend on $k$.

\cvb
The $C^2$- estimates in \rs{5} remain valid, too, since the second derivatives of
$ \tilde f$, $ \tilde f^j_i$, that occur on the right-hand side of \re{3.14}, can be
expressed as---we only consider the covariant form $ \tilde f_{ii}$, no summation
over $i$---

\begin{equation}
\begin{aligned}
- \tilde f_{ii}= - \tilde f_{\al\bet} x^\al_i x^\bet_i &-2 \tilde f_{\al \tilde \n^\bet}
x^\al_i
\tilde
\n^\bet_i \\[\cma]
&- \tilde f_\al\n^\al h_{ii}- \tilde f_{ \tilde \n^\al \tilde \n^\bet}
\tilde
\n^\al_i \tilde \n^\bet_i -\tilde f_{\tilde\n^\al}\tilde\n^\al_{i;i},
\end{aligned}
\end{equation}

\cvm
\nd where

\begin{equation}
\tilde \n^\al_i=\tht \tilde v^{-1} \n^\al_i+\dot\tht \tilde v_i \tilde v^{-1} \n^\al -\tht
\tilde v^{-2} \tilde v_i\n^\al,
\end{equation}

\begin{equation}
\begin{aligned}
\tilde \n^\al_{i;i}&=2\dot\tht \tilde v_i \tilde v^{-1}\n^\al_i - 2\tht \tilde v^{-2}
\tilde v_i\n^\al_i + \tht \tilde v^{-1} \n^\al_{i;i}\\[\cma]
&\hp{=} +\Ddot\tht \tilde v_i \tilde v_i \tilde v^{-1}\n^\al - 2\dot\tht \tilde v_i
\tilde v_i \tilde v^{-2}\n^\al+\dot\tht \tilde v_{ii}\tilde v^{-1}\nu^\al\\[\cma]
&\hp{=}+ 2\tht \tilde v^{-3} \tilde v_i \tilde v_i\n^\al-\tht \tilde v^{-2} \tilde v_{ii}
\n^\al,
\end{aligned}
\end{equation}

\nd and

\begin{equation}
\begin{aligned}
\tilde v_{ii}=\h_{\al\bet\ga}x^\bet_i x^\ga_i \n^\al +\h_{\al\bet} \n^\al \n^\bet h_{ii}+2
\h_{\al\bet} x^\bet_i \n^\al_i +\h_\al \n^\al_{i;i}.
\end{aligned}
\end{equation}

\cvm
Hence, the result of \rp{5.1} is still valid since no additional bad terms occur in
inequality \re{5.9} as one easily checks, and since, furthermore, we also have

\begin{equation}
\tilde f\le F_\e
\end{equation}

\cvm
\nd during the evolution, for the modified version of \re{3.10} now has the form

\begin{equation}
\begin{aligned}
{( F_\e- \tilde f)}^\prime&- F^{ij}_\e( F_\e- \tilde f)_{ij}\\[\cma]
&= - F^{ij}_\e h_{ik}h_j^k ( F_\e- \tilde f)
 - \tilde f_\al\n^\al ( F_\e- \tilde f)\\[\cma]
&\hp{=\;} - \tilde f_{ \tilde \n^\ga}\n^\ga [\dot\tht \tilde v^{-1}-\tht \tilde v^{-2}]
\h_{\al\bet} \n^\al \n^\bet (F_\e- \tilde f)\\[\cma]
&\hp{=\; }  -[\dot\tht \tilde v^{-1}-\tht \tilde v^{-2}] \tilde f_{ \tilde
\n^\bet}\n^\bet \h_\al x^\al_i (F_\e- \tilde f)_jg^{ij}\\[\cma]
&\hp{+\; }-\tht \tilde v^{-1} \tilde f_{ \tilde \n^\al} x^\al_i (F_\e - \tilde f)_j
g^{ij}\\[\cma]
&\hp{=\; }   - F^{ij}_\e\riema
\al\bet\ga\de\n^\al x_i^\bet
\n^\ga x_j^\de ( F_\e-\tilde  f).
\end{aligned}
\end{equation}

\cvm
Here, we used the relation

\begin{equation}
\dot{\tilde v}=\h_{\al\bet} \dot x^\bet\n^\al+\h_\al\dot\n^\al,
\end{equation}

\cvm
\nd which follows immediately from \re{4.8}, together with \re{3.8} and \re{3.7}.

\cvm
The conclusions of \rs{6} are therefore applicable leading to  a solution of
equation
\re{7.5}.

\cvb
\section{$C^1$- estimates for the stationary approximations}\las{8}

\cvb
Consider the solutions $M_\e=\graph u_\e$ of equation \re{7.5}, which at the
moment not only depend on $\e$ but also on $k$, the parameter of the cut-off
function $\tht$, cf. \re{7.3}. We shall prove that the hypersurfaces $M_\e$ are
uniformly space-like independent of $\e$ and $k$, or, equivalently, that there
exists a  constant $m_1$ such that

\begin{equation}\lae{8.1}
\tilde v=(1-\abs{Du_\e}^2)^{-\frac12}\le m_1\qq \A\; \e, k,
\end{equation}

\cvm
\nd where the parameter $\e$ is supposed to be small and $k$ to be large, so
that the barrier condition is satisfied.

\cvb
\bl
Let $u_\e$ be a solution of \re{7.5}, then, the estimate \re{8.1} is valid
uniformly in $\e$ and $k$. Hence, 
$M_\e=\graph u_\e$ is a solution of equation \re{7.2}, if we choose $k\ge 2
m_1$.
\el

\cvb
\bp
For arbitrary but fixed values of $\e$ and $k$, let us introduce the notation $
\tilde F$ for $F_\e$, where from now on through the rest of the article

\begin{equation}
F=H_2,
\end{equation}

\cvm
\nd and where $f=f(x,\n)$ satisfies 

\begin{align}
0<c_1&\le f(x,\n),\lae{8.3}\\[\cma]
\nnorm{f_\bet(x,\n)}&\le c_2(1+ \nnorm\n^2),\lae{8.4}\\
\intertext{and}
\nnorm{f_{\n^\bet}(x,\n)}&\le c_3 (1+\nnorm\n),\lae{8.5}
\end{align}

\cvm
\nd for all $x\in\bar\Om$ and all past directed time-like vectors $\n\in T_x(\Om)$.

\cvm
Thus, $F$ is homogeneous of degree 2, and we recall that

\begin{equation}\lae{8.6}
F^{ij}=H g^{ij}-h^{ij},
\end{equation}

\nd and
\begin{equation}\lae{8.7}
\tilde F^{ij}=F^{ij}+\e (n-1) (1+\e n) H g^{ij},
\end{equation}

\cvm
\nd where $ \tilde F^{ij}$ is evaluated at $h_{ij}$ and $F^{ij}$ at $(h_{ij}+ \e H
g_{ij})$.

 \cvm
We also drop the index $\e$, writing $u$ for $u_\e$ and $M$ for $M_\e$, i.e.
$M$ solves the equation

\begin{equation}
\fv{\tilde F}M=f(x, \tilde \n).
\end{equation}

\cvb
The $C^1$- estimate will follow the arguments in the proof of \rp{4.8}, where at
one point we shall introduce an additional observation especially suitable for the
curvature function $F=H_2$.

\cvb
\br\lar{8.2}
The former parabolic equations and inequalities, \re{4.10}, \re{4.39b}, and
\re{4.47b} can now be read as elliptic equations resp. inequalities
by simply assuming that the terms involved are time independent. Though, to
be absolutely precise, one has to observe that the present curvature function is
homogeneous of degree 2, which means that, whenever the term $F$---not
derivatives of $F$---occurs explicitly in the  equations or inequalities
just mentioned, it has to be replaced by $2F$ because it was obtained as a result
of Euler's formula for homogeneous functions of degree $d_0$

\begin{equation}
d_0 F=F^{ij} h_{ij}.
\end{equation}

\cvm
We mention it as a matter of fact only, since it doesn't affect the estimates at all.
\er

\cvb
However, we have to be aware that $f$ now depends on $ \tilde \n$ instead of
$\n$, i.e. the elliptic version of inequality \re{4.47b} now takes the form

\begin{equation}\lae{8.10}
\begin{aligned}
 -\tilde F^{ij} w_{ij}&\le -\de\msp[2] \tilde F^{ij} h^k_i h_{kj} \tilde v \f +
c[\de^{-1}+\abs\lam \m e^{\lam u}] \tilde F^{ij} g_{ij} \tilde v^3 \f\\[\cma]
&\hp{\le} +[\frac1{1-2\de}-1] \tilde F^{ij} u_i u_j \m^2\lam^2 e^{2\lam u} \tilde v
\f \\[\cma]
&\hp{\le}-\tilde F^{ij} u_i u_j \m\lam^2 e^{\lam u} \tilde v \f \\[\cma]
&\hp{\le} +2f
[-\h_{\al\bet}\n^\al \n^\bet+e^{-\psi}
\m
\lam e^{\lam u} \tilde v^2] \f + c\msp[2] \nnorm{f_\bet} \tilde v^2 \f \\[\cma]
&\hp{\le}  +f_{ \tilde \n^\bet}[\tht \tilde v ^{-1}\n^\bet_i +\dot\tht \tilde v _i \tilde
v^{-1} \n^\bet -\tht \tilde v^{-2} \tilde v_i \n^\bet] u^i e^\psi 
\f.
\end{aligned}
\end{equation}

\cvm
Here, we used the notation $\de=\de(x)$ for the small parameter in the Schwarz
inequality instead of
$\e$, which has a different meaning in the present context, $w$ is defined as in
\re{4.44}, where the parameters
$\m,
\lam$ should  satisfy the conditions in \re{4.45}, and $u$ is supposed to be less than
$-1$.

\cvm
We claim that $w$ is uniformly bounded  provided $\m$ and $\lam$ are chosen
appropriately. Following the arguments  in \rs{4}, we shall use the maximum
principle and consider a point $x_0\in M$, where

\begin{equation}
w(x_0)=\sup_M w.
\end{equation}

\cvm
 As before, we   choose
$\de=e^{-\lam u}$. But the further conclusions are no longer valid, since we have a
really bad term on the right-hand side of \re{8.10} that is of the order $ \tilde
v^4$ due to the assumption \re{8.4}.

\cvm
The only possible good term which can  balance it, is

\begin{equation}\lae{8.12}
 -\de\msp[2] \tilde  F^{ij} h^k_i h_{kj} \tilde v \f.
\end{equation}

\cvm
To exploit this term we use the fact that $Dw(x_0)=0$, or, equivalently

\begin{equation}\lae{8.13}
\begin{aligned}
- \tilde v_i&= \m \lam e^{\lam u} \tilde v u_i\\[\cma]
&= e^\psi h^k_iu_k -\h_{\al\bet}\n^\al x^\bet_i,
\end{aligned}
\end{equation}

\cvm
\nd where the second equation follows from  \re{4.8} and the definition of
the covariant vectorfield
$\h=e^\psi (-1,0,\dotsc,0)$.

\cvm
Next, we choose a coordinate system $(\x^i)$ such that in the critical point

\begin{equation}\lae{8.14}
g_{ij}=\de_{ij}\qq\tup{and}\qq h^k_i=\ka_i \de^k_i,
\end{equation}

\cvm
\nd and the labelling of the principal curvatures corresponds to

\begin{equation}\lae{8.15}
\ka_1\le \ka_2\le \dotsb \le \ka_n.
\end{equation}

\cvm
Then, we deduce from \re{8.13}

\begin{equation}
e^\psi \ka_i u_i=\m\lam e^{\lam u} \tilde v u_i + \h_{\al\bet}\n^\al x^\bet_i.
\end{equation}

\cvm
Assume that  $ \tilde v(x_0)\ge 2$, and let $i=i_0$ be an index such that

\begin{equation}\lae{8.17}
\abs{u_{i_0}}^2\ge \frac1n\norm{Du}^2.
\end{equation}

\cvm
Setting $(e^i)=\pde{}{\x^{i_0}}$ and assuming without loss of generality that
$0\le u_ie^i$ in $x_0$ we infer from \rl{2.6}

\begin{equation}
\begin{aligned}
e^\psi \ka_{i_0} u_ie^i&=\m\lam e^{\lam u} \tilde v u_i e^i+\h_{\al\bet}\n^\al
x^\bet_ie^i \\[\cma]
&\le \m\lam e^{\lam u} \tilde v u_ie^i +c \tilde v^2,
\end{aligned}
\end{equation}

\cvm
\nd and we deduce further in view of \re{2.12b} and \re{8.17} that 

\begin{equation}\lae{8.19}
\begin{aligned}
\ka_{i_0}\le [\m\lam e^{\lam u}+ c] \tilde v e^{-\psi} 
\le \frac12\m\lam e^{\lam u} \tilde  v e^{-\psi},
\end{aligned}
\end{equation}

\cvm
\nd if $\abs\lam$ is sufficiently large,
 i.e. $\ka_{i_0}$ is \tit{negative} and of the same order as $ \tilde v$.

\cvm
The Weingarten equation and \rl{2.6} yield

\begin{equation}
\nnorm{\n^\bet_iu^i}=\nnorm{h^k_iu^i x^\bet_k}\le c \tilde v [h^k_iu^i
h_{kl}u^l]^{\frac12},
\end{equation}

\cvm
\nd and therefore, we infer from \re{8.13}

\begin{equation}
\abs{ \tilde v_i u^i}+\nnorm{\n^\bet_i u^i}\le c \m\abs\lam e^{\lam u} \tilde v^3
\end{equation}

\cvm
\nd in critical points of $w$, and hence, that in those points, the term involving
$f_{\tilde \n^\bet}$ on the right-hand side of inequality \re{8.10} can be estimated
from above by

\begin{equation}\lae{8.21}
\abs{f_{ \tilde \n^\bet}[\tht \tilde v ^{-1}\n^\bet_i +\dot\tht \tilde v _i \tilde
v^{-1} \n^\bet -\tht \tilde v^{-2} \tilde v_i \n^\bet] u^i e^\psi 
\f}\le c \m\abs\lam e^{\lam u} \tilde v^4\f.
\end{equation}

\cvm
Next, let us estimate the crucial term in \re{8.12}. Using \re{8.7}, the particular
coordinate system \re{8.14}, as well as the inequalities \re{8.15}, together with
the fact that $\ka_{i_0}$ is negative, we conclude

\begin{equation}\lae{8.22}
\begin{aligned}
- \tilde F^{ij}h^k_ih_{kj}&\le - F^{ij}h^k_i h_{kj}\le -\sum_{i=1}^{i_0} F^i_i
\ka^2_i\\[\cma]
&\le -\sum_{i=1}^{i_0} F^i_i \ka^2_{i_0},
\end{aligned}
\end{equation}

\cvm
\nd where we recall that the argument of $F^i_i$ is the $n$- tupel with
components

\begin{equation}\lae{8.23}
\tilde \ka_j=\ka_j+\e H
\end{equation}

\cvm
\nd and observe that in the present coordinate system 

\begin{equation}
F^i_i=\pde F{ \tilde \ka_i}.
\end{equation}

\cvm
Let $\hat F=\log F$, then, $ \hat F$ is concave, and therefore, we have in
view of \re{8.15}

\begin{equation}
\hat F^1_1\ge \hat F^2_2\ge\dotsb\ge \hat F^n_n,
\end{equation}

\cvm
\nd cf. \ci[Lemma 2]{eh}, or equivalently,

\begin{equation}
 F^1_1\ge F^2_2\ge\dotsb\ge  F^n_n.
\end{equation}

\cvm
Hence, we conclude

\begin{equation}\lae{8.27}
\begin{aligned}
-\sum_{i=1}^{i_0}F^i_i&\le -F^1_1\le -\frac1n\sum_{i=1}^nF^i_i\\[\cma]
&= -\frac1n(n-1)[H+\e n H]\le -\frac{n-1}n H,
\end{aligned}
\end{equation}

\cvm
\nd where we also used \re{8.6} and \re{8.23}.

\cvm
Combining \re{8.19}, \re{8.22}, \re{8.27}, and the estimate \re{1.15}, we
deduce further

\begin{equation}
\begin{aligned}
- \tilde F^{ij} h^k_i h_{kj}&\le  -\frac{n-1}n H \ka_{i_0}^2\\[\cma]
&\le -\frac{n-1}{2n} c_n \abs{\ka_{i_0}}^3-\frac{n-1}{2n} H \ka_{i_0}^2\\[\cma]
&\le -a_0 \m^3\abs\lam^3 e^{3\lam u} \tilde v^3-a_1 H \m^2\lam^2 e^{2\lam u}
\tilde v^2
\end{aligned}
\end{equation}

\cvm
\nd with some positive constants $a_0=a_0(n,\Om)$ and $a_1=a_1(n,\Om)$.

\cvm
Inserting this estimate, and the estimate in \re{8.21} in the elliptic
inequality \re{8.10}, with $\de=e^{-\lam u}$, we finally obtain

\begin{equation}
\begin{aligned}
 -\tilde F^{ij} w_{ij}&\le -a_0\m^3\abs\lam^3 e^{2\lam u} \tilde v^4 \f -a_1 H
\m^2\lam^2 e^{\lam u}
\tilde v^3\f\\[\cma]
&\hp{\le}
+\frac2{1-2\de}\msp[2]\tilde F^{ij} u_i u_j \m^2\lam^2 e^{\lam u} \tilde v
\f +c [1+\abs\lam \m] H e^{\lam u} \tilde v^3\f\\[\cma]
&\hp{\le}-\tilde F^{ij} u_i u_j \m\lam^2 e^{\lam u} \tilde v \f +c[c_2+c_3 \m\abs\lam
e^{\lam u}] \tilde v^4
\f\\[\cma]
&\hp{\le} +2f
[c +e^{-\psi}
\m
\lam e^{\lam u} ] \tilde v^2\f.  
\end{aligned}
\end{equation}

\cvm
Choosing, now, $\m=\frac14$ and $\abs\lam$ large, the right-hand side of the
preceding inequality is negative, contradicting the maximum principle, i.e. the
maximum of
$w$ cannot occur at point where $ \tilde v\ge 2$. Thus, the desired uniform
estimate for $w$ and, hence, $ \tilde v$ is proved.
\ep

\cvb
Let us close this section with an interesting observation that is an immediate
consequence of the preceding proof, we have especially  \re{8.22} and
the first line of inequality \re{8.27} in mind,

\cvm
\bl\lal{8.3}
Let $F\in C^2(\C)$ be a positive symmetric curvature function such that the
partial derivatives $F_i$ are positive and $ \hat F=\log F$ is concave. Suppose
$F$ is evaluated at a point $(\ka_i)$, and assume that $\ka_{i_0}$ is a
component  that is either negative or the smallest component of that particular
$n$- tupel, then

\begin{equation}\lae{8.30}
\sum_{i=1}^n F_i\ka^2_i\ge \frac1n \sum_{i=1}^n F_i \msp[2] \ka^2_{i_0}.
\end{equation}
\el

\cvb
\section{$C^2$- estimates for the stationary approximations}\las{9}

\cvb
We want to prove uniform $C^2$- estimates for admissible solutions $M$ of

\begin{equation}\lae{9.1}
\fv{ \tilde F}{M}=f(x,\n),
\end{equation}

\cvm
\nd where we use the notations and conventions of the preceding section.

The starting point is  an elliptic equation for the second fundamental form.

\bl
The tensor $h^j_i$ satisfies the elliptic equation

\begin{equation}\lae{9.2}
\begin{aligned}
- &\tilde F^{kl}h_{i;kl}^j\\
&=-  \tilde F^{kl}h_{rk}h_l^rh_i^j
+2f h_i^kh_k^j\\
&\hp{+}- f_{\al\bet} x_i^\al x_k^\bet g^{kj}-  f_\al\n^\al
h_i^j-f_{\al\n^\bet}(x^\al_i x^\bet_kh^{kj}+x^\al_l
x^\bet_k h^{k}_i\, g^{lj})\\ &\hp{=}
-f_{\n^\al\n^\bet}x^\al_lx^\bet_kh^k_ih^{lj}-f_{\n^\bet} x^\bet_k h^k_{i;l}\,g^{lj} 
-f_{\n^\al}\n^\al h^k_i h^j_k\\ &\hp{=}+
\tilde F^{kl,rs}h_{kl;i}h_{rs;}^{\hphantom{rs;}j}+2 \tilde F^{kl}\riema
\al\bet\ga\de x_m^\al x_i ^\bet x_k^\ga x_r^\de h_l^m g^{rj}\\
&\hp{=}- \tilde F^{kl}\riema \al\bet\ga\de x_m^\al x_k ^\bet x_r^\ga x_l^\de
h_i^m g^{rj}- \tilde F^{kl}\riema \al\bet\ga\de x_m^\al x_k ^\bet x_i^\ga x_l ^\de
h^{mj} \\ &\hp{=}- \tilde F^{kl}\riema \al\bet\ga\de\n^\al x_k^\bet\n^\ga x_l^\de
h_i^j+ 2f
\riema \al\bet\ga\de\n^\al x_i^\bet\n^\ga x_m^\de g^{mj}\\
&\hp{=}+ \tilde F^{kl}\bar R_{\al\bet\ga\de;\e}\{\n^\al x_k^\bet x_l^\ga x_i^\de
x_m^\e g^{mj}+\n^\al x_i^\bet x_k^\ga x_m^\de x_l^\e g^{mj}\}.
\end{aligned}
\end{equation}
\el

\bp
The elliptic equation can be immediately derived from the corresponding
parabolic equation \re{3.14} by dropping the time-derivative,  replacing $F$
by $ \tilde F$, and observing that the present curvature function is
homogeneous of degree 2, cf. \rr{8.2}.
\ep

\cvb
Contracting over the indices $(i,j)$ in \re{9.2} we obtain a differential equation
for $H$

\begin{equation}\lae{9.3}
\begin{aligned}
- &\tilde F^{kl}H_{kl}\\[\cma]
&=-  \tilde F^{kl}h_{rk}h_l^rH
+2f \abs A^2\\[\cma]
&\hp{+}- f_{\al\bet} x_i^\al x_k^\bet g^{ki}-  f_\al\n^\al
H-2f_{\al\n^\bet}x^\al_i x^\bet_kh^{ki}\\[\cma]
&\hp{=}
-f_{\n^\al\n^\bet}x^\al_lx^\bet_kh^k_ih^{li}-f_{\n^\bet}(H^k x^\bet_k +\abs
A^2\n^\bet)\\[\cma]
&\hp{=}-\bar R_{\al\bet}\n^\al x^\bet_kg^{kl}x^\ga_l f_{\n^\ga} 
+
\tilde F^{kl,rs}h_{kl;i}h_{rs;}^{\hphantom{rs;}i}\\[\cma]
&\hp{=}+2 \tilde F^{kl}\riema
\al\bet\ga\de x_m^\al x_i ^\bet x_k^\ga x_r^\de h_l^m g^{ri} - 2\tilde F^{kl}
\riema \al\bet\ga\de x_m^\al x_k ^\bet x_i^\ga x_l ^\de h^{mi}\\[\cma]
&\hp{=}- \tilde F^{kl}\riema \al\bet\ga\de\n^\al x_k^\bet\n^\ga x_l^\de
H+ 2f
\bar R_{\al\bet}\n^\al\n^\bet\\[\cma]
&\hp{=}+ \tilde F^{kl}\bar R_{\al\bet\ga\de;\e}\{\n^\al x_k^\bet x_l^\ga x_i^\de
x_m^\e g^{mi}+\n^\al x_i^\bet x_k^\ga x_m^\de x_l^\e g^{mi}\},
\end{aligned}
\end{equation}

\cvm
\nd where we also used the symmetry properties of the Riemann curvature tensor
and  the Codazzi equations at one point.

\cvm
Next, let us improve the estimate in \rp{5.1}.

\bl\lal{9.2}
Let $M=\graph u$ be an admissible solution of \re{9.1}, then the principal
curvatures of $M$ satisfy the estimate

\begin{equation}\lae{9.4}
\e\abs A^2\le \const,
\end{equation}

\cvm
\nd where the constant depends on $\nnorm{Df}, \nnorm{D^2f}$,  the
constant $c_1$ in \re{8.3}, and on known estimates of the $C^0$ and $C^1$-
norm of $u$.
\el

\bp
We argue as in the proof of \rp{5.1} and define

\begin{equation}
\f=\sup\set{{h_{ij}\h^i\h^j}}{{\norm\h=1}}.
\end{equation}

\cvm
Let $x_0\in M$ be a point, where $\f$ achieves its maximum, and assume
without loss of generality that, after having introduced normal Riemannian
coordinates around $x_0$, we may write $\f=h^n_n$, cf. the corresponding
arguments in the proof of \rp{5.1}.

\cvm
Applying the maximum principle in $x_0$, we deduce from \re{9.2} the following
inequality

\begin{equation}
\begin{aligned}
0&\le -\e (n-1) H \abs A^2 h^n_n+2 f\abs{h^n_n}^2+ \tilde F^{kl,rs}h_{kl;n}
h_{rs;}^{\hp{rs;}n}\\[\cma] &\hp{\le\;} +c(1+\nnorm{Df}+\nnorm{D^2f})(1+\abs
A^2)+c ( \tilde F^{ij}g_{ij}+f),
\end{aligned}
\end{equation}

\cvm
\nd where we also used \re{8.7}, the Weingarten and Codazzi equations, and the
fact that the pair $(x,\n)$ stays in a compact subset of $\bar\Om\times
C_-(\bar\Om)$, where $C_-(\bar\Om)$ stands for the set of past directed time-like
vectorfields in
$\bar\Om$.

\cvm
Furthermore, we know that

\begin{equation}
\tilde F^{ij}g_{ij}\le c H,
\end{equation}
and
\begin{equation}
\begin{aligned}
\tilde F^{kl,rs}h_{kl;n}h_{rs;}^{\hp{rs;}n} &\le \tilde F^{-1} \tilde F_{;n} \tilde
F^{\hp{;}n}_;=f^{-1} f_nf^n\\[\cma]
&\le c\msp[2] c_1^{-1} (1+\abs A^2),
\end{aligned}
\end{equation}

\cvm
\nd since $\log \tilde F$ is concave, cf. \rl{1.4}.

\cvm
Thus,  we conclude that in $x_0$ the following inequality is valid

\begin{equation}
0\le -\e \abs A^4 +c(1+\abs A^2)
\end{equation}

\cvm
\nd with a known constant $c$, and the lemma is proved.
\ep

\cvb
The estimate \re{9.4} will play an important role in the final a priori estimate.

\cvb
\bl\lal{9.3}
Let $F=H_2$, $M=M^n$ a Riemannian manifold with metric $g_{ij}$,  $h_{ij}$ a
symmetric tensor field on $M$ the eigenvalues of which belong to $\C_2$, and
$p\in M$ an arbitrary point. Choose local coordinates around $p$ such that the
relations \re{8.14} and \re{8.15} are satisfied. Then, we have for $1\le j\le n$

\begin{align}
\sum_{i\ne j} \ka^2_i + 2F&=\abs{F^j_j}^2+2F^j_j\ka_j,\lae{9.10}\\[\cma]
\sum_{i\ne n}\ka^2_i+2F&\le c F^n_n\ka_n,\lae{9.11}
\end{align}
and
\begin{equation}\lae{9.12}
\sum_{i\ne
j}\absb{\frac{F_{;j}}{F^j_j}+\frac{F_{;j}}H(1-\frac{F^i_i}{F^j_j})}^2\le
c\abs{F_{;j}}^2 F^{-1},
\end{equation}

\cvm
\nd with $c=c(n)$, $F$ is evaluated at $h_{ij}$, and where we point out that the
summation convention is not used.
\el

\cvb
\bp
Throughout the proof we shall use the ambivalent meaning of $F$ as a function
depending on $\ka_i$ or on $h_{ij}$ switching freely from one viewpoint  to the
other.

\cvm (i) From the definition of $F$

\begin{equation}
F=\frac12(H^2-\abs A^2),
\end{equation}

\cvm
\nd and \re{8.6} we conclude

\begin{equation}
\begin{aligned}
2F&=(F^j_j+\ka_j)^2-\sum_i\ka^2_i\\[\cma]
&=\abs{F^j_j}^2+2F^j_j\ka_j-\sum_{i\ne j}\ka^2_i,
\end{aligned}
\end{equation}

\cvm
\nd which proves \re{9.10}.

\cvm
(ii) If $j=n$, and thus $\ka_n$ the largest eigenvalue, then, we derive from
\re{8.6}

\begin{equation}
F^n_n\le H\le n\ka_n,
\end{equation}

\cvm
\nd and \re{9.11} follows at once from \re{9.10}.

\cvm
(iii) A simple algebraic transformation yields

\begin{equation}\lae{9.16}
\begin{aligned}
\frac{F_{;j}}{F^j_j}+\frac{F_{;j}}H(1-\frac{F^i_i}{F^j_j})&=
\frac{F_{;j}}{HF^j_j}(H-\ka_j+\ka_i)\\[\cma]
&=
\frac{F_{;j}}{HF^j_j}(\sum_{k\ne j}\ka_k+\ka_i),
\end{aligned}
\end{equation}

\cvm
\nd and hence,

\begin{equation}\lae{9.17}
\begin{aligned}
\sum_{i\ne
j}\absb{\frac{F_{;j}}{F^j_j}+\frac{F_{;j}}H(1-\frac{F^i_i}{F^j_j})}^2\le c\msp[2]
\frac{\abs{F_{;j}}^2}{H^2\abs{F^j_j}^2}\sum_{i\ne j}\ka^2_i.
\end{aligned}
\end{equation}

\cvm
We, now, treat the cases $j=n$ and $j\ne n$ separately. 

\cvm
If $j=n$, we apply
\re{9.11} and \re{1.17} and conclude

\begin{equation}
\begin{aligned}
\sum_{i\ne
j}\absb{\frac{F_{;j}}{F^j_j}+\frac{F_{;j}}H(1-\frac{F^i_i}{F^j_j})}^2\le c\msp[2]
\frac{\abs{F_{;j}}^2}{HF^j_j}\le c\msp[2]
\frac{\abs{F_{;j}}^2}{F}.
\end{aligned}
\end{equation}

\cvm
If $j\ne n$, we deduce from \re{9.10}

\begin{equation}
\ka^2_n\le 6\abs{F^j_j}^2,
\end{equation}

\cvm
\nd and deduce further

\begin{equation}
\sum_{i\ne j} \ka^2_i \le 8\abs{F^j_j}^2,
\end{equation}

\cvm
\nd where apparently we only had to worry about the case $0\le \ka_j$.

\cvm
Thus, the right-hand side of \re{9.17} is estimated from above by

\begin{equation}
c\msp[2]\frac{\abs{F_{;j}}^2}{H^2},
\end{equation}

\cvm
\nd which in turn is less than

\begin{equation}
c\msp[2]\frac{\abs{F_{;j}}^2}F.
\end{equation}
\ep

\cvb
\bc
Let $M$ be an admissible solution of \re{9.1} and $p\in M$ arbitrary. Choose
local coordinates around $p$ such that the relations \re{8.14} and \re{8.15} are
valid. Then, for any $1\le j\le n$, the following inequality is valid in $p$

\begin{equation}\lae{9.23}
\sum_{i\ne
j}\absb{\frac{f_j}{\tilde F^j_j}+\frac{f_j}{ \widetilde
H}(1-\frac{\tilde F^i_i}{\tilde F^j_j})}^2\le c\abs{f_j}^2 f^{-1},
\end{equation}

\cvm
\nd where we use the notation $ \widetilde H=(1+\e n)H$, and do not apply the
summation convention.
\ec

\cvb
\bp
Let us recall the relation \re{8.7}, which we can also express in the form

\begin{equation}\lae{9.24}
\tilde F^{ij}= \widetilde H g^{ij}- \tilde h^{ij}+\e F^{rs}g_{rs} g^{ij},
\end{equation}

\cvm
\nd where 

\begin{equation}\lae{9.25}
\tilde h_{ij}=h_{ij}+\e Hg_{ij}.
\end{equation}

\cvm
Consider each summand in \re{9.23} separately. We have

\begin{equation}
\begin{aligned}
\absb{\frac{f_j}{\tilde F^j_j}+\frac{f_j}{ \widetilde
H}(1-\frac{\tilde F^i_i}{\tilde F^j_j})}^2&=\frac{f_j^2}{ \widetilde H ^2
\abs{\tilde F^j_j}^2}\abs{ \widetilde H + \tilde F^j_j- \tilde F^i_i}^2\\[\cma]
&\le \frac{f_j^2}{ \widetilde H ^2
\abs{ F^j_j}^2}\absb{\sum_{k\ne j} \tilde \ka_k+ \tilde \ka_i}^2,
\end{aligned}
\end{equation}

\cvm
\nd in view of \re{9.24}, i.e. we are exactly in the same situation as in the proof of
\rl{9.3} after the equation \re{9.16} with the following modifications: we replace
$h_{ij},
\ka_i$ and $H$ by $ \tilde h_{ij}, \tilde \ka_i$ resp. $ \widetilde H$ and observe
that
$F( \tilde h_{ij})=f$.
\ep

\cvb
\bl\lal{9.5}
Let $M$ be an admissible solution of equation \re{9.1}, then, the estimate

\begin{multline}\lae{9.27}
\tilde F^{ij,kl}h_{ij;r}h^{\hp{kl;}r}_{kl;}  +H^{-1} \tilde F^{ij} H_i H_j
\le c f^{-1}
\norm{Df}^2 + c \msp \e \norm{DH}^2+c
\end{multline}

\cvm
\nd is valid in every point, where the smallest principal curvature $\ka_1$ satisfies

\begin{equation}\lae{9.28}
\max(-\ka_1,0)\le \frac1{2(n-1)} H\equiv \e_1 H.
\end{equation}
\el

\cvm
\bp
The proof is a modification of a similar result in \ci[Sections 6.1.4 and 6.1.5]{pb}
or \ci[Section 5.1.1]{pb2}.

\cvm
It follows immediately from the definition of $F$ that

\begin{equation}\lae{9.29}
F^{ij,kl}= g^{ij} g^{kl} -\frac12(g^{ik}g^{jl}+g^{il}g^{jk})
\end{equation}
\nd and
\begin{equation}\lae{9.30}
\begin{aligned}
\tilde F^{ij,kl}h_{ij;r} h_{kl;s} g^{rs}&= F^{ij,kl}h_{ij;r} h_{kl;s} g^{rs}\\[\cma]
&\hp{=\; }+ \e[2(n-1) +\e(n-1)n] \norm{DH}^2.
\end{aligned}
\end{equation}

\cvm
In a fixed point $p\in M$ introduce normal Riemannian coordinates such that the
relations \re{8.14} and \re{8.15} are valid, and define the matrix $(a^{kl})$
through

\begin{equation}
a^{kl}=\begin{cases}
 1, &k\ne l,\\
0 , &k=l.
\end{cases}
\end{equation}

\cvm
We also set

\begin{equation}
h_{ijk}=h_{ij;k}.
\end{equation}

\cvm
Then, we conclude from \re{9.29}

\begin{equation}\lae{9.33}
\begin{aligned}
F^{ij,kl}h_{ijr} h_{kls} g^{rs}&= \norm{DH}^2 - h_{ijk}h^{ijk}\\[\cma]
&= \sum_{i,k,l} (h_{kki} h_{lli}-h_{kli}h_{kli})\\[\cma]
&= \sum_i a^{kl}(h_{kki} h_{lli}- h^2_{kli})\\[\cma]
&= \sum_i a^{kl} h_{kki} h_{lli}- \sum_l\sum_{k,i} a^{kl} h^2_{kli}.
\end{aligned}
\end{equation}

\cvm
In the last summand let us interchange the roles  of $i$ and $l$ to obtain

\begin{equation}
-\sum_l\sum_{k,i} a^{kl} h^2_{kli}=-\sum_i\sum_{k,l} a^{ki} h^2_{kil}.
\end{equation}

\cvm
The Codazzi equations yield

\begin{equation}
h_{kil}= h_{kli} + c_{kil},
\end{equation}

\cvm
\nd where $(c_{kil})$ is a uniformly bounded tensor in $\bar\Om$

\begin{equation}
\nnorm{c_{kil}}\le \const,
\end{equation}

\cvm
\nd and hence,

\begin{equation}
h^2_{kil}=h^2_{kli}+2 h_{kli} c_{kil} + c^2_{kil}
\end{equation}
and
\begin{equation}\lae{9.38}
\begin{aligned}
-\sum_i\sum_{k,l} a^{ki} h^2_{kil}&\le -\sum_i\sum_{k,l} a^{ki}
h^2_{kli}-2\sum_{i,k,l} a^{ki} h_{kli} c_{kil} \\[\cma]
&\le - 2\sum_i\sum_k a^{ki} h^2_{kki} - \sum_{[i,j,k]} h^2_{ijk}\\[\cma]
&\hp{\le\; }-2\sum_{i,k,l} a^{ki} h_{kli} c_{kil}-2\sum_{i,k}a^{ki}h_{iik}c_{iki},
\end{aligned}
\end{equation}

\cvm
\nd where $\sum_{[i,j,k]}$ means that the summation is carried out over those
triples $(i,j,k)$ where all three indices are different from each other.

\cvm
The first linear term in the last inequality can be estimated from above by

\begin{equation}
\begin{aligned}
-2\sum_{i,k,l} a^{ki} h_{kli} c_{kil}&= -2 \sum_i\sum_k a^{ki} h_{kki} c_{kik}-
2 \sum_i\sum_{k\ne l} a^{ki} h_{kli} c_{kil}\\[\cma]
&\le \frac{\de}2 \sum_i\sum_k a^{ki} h^2_{kki}+\frac{\de}2 \sum_{[i,j,k]} h^2_{ijk}
+c\de^{-1}
\end{aligned}
\end{equation}

\cvm
\nd for any $\de>0$, and the second term similarly. Thus,  we deduce from
\re{9.33} and
\re{9.38}

\begin{equation}\lae{9.40}
\begin{aligned}
F^{ij,kl}h_{ijr} h_{kls} g^{rs}&\le -2(1-\frac\de2) \sum_i\sum_k a^{ki}
h^2_{kki}\\[\cma]
&\hp{\le\; }+\sum_i a^{kl} h_{kki}h_{lli} + c\de^{-1}
\end{aligned}
\end{equation}

\cvm
\nd for any $0<\de\le 1$.

\cvm
Next, let us consider the second term on the left-hand side of \re{9.27}; we have

\begin{equation}
\begin{aligned}
H^{-1} \tilde F^{ij} H_i H_j&= \widetilde H^{-1} \tilde F^{ij} H_i H_j +\frac{\e
n}{1+\e n} H^{-1} \tilde F^{ij} H_i H_j\\[\cma]
&\le \widetilde H^{-1} \sum_i \tilde F^{ii} (\sum_kh_{kki})^2 +\e c
\norm{DH}^2,
\end{aligned}
\end{equation}

\cvm
\nd where $c=c(n)$, in view of \re{8.7},  \re{9.25}, and \re{1.15}. 

\cvm
Combining \re{9.40} and the preceding estimate, we conclude

\begin{equation}\lae{9.42}
\begin{aligned}
F^{ij,kl}&h_{ij;r}h_{kl;s} g^{rs} +H^{-1} \tilde F^{ij} H_i H_j\\[\cma]
&\le 
-2(1-\frac\de2) \sum_i\sum_k a^{ki}
h^2_{kki}
+\sum_i a^{kl} h_{kki}h_{lli} + c\de^{-1}\\[\cma]
&\hp{\le\; }+\widetilde H^{-1} \sum_i \tilde F^{ii} (\sum_kh_{kki})^2 +\e c
\norm{DH}^2.
\end{aligned}
\end{equation}

\cvm
For each index $i$, let us estimate the corresponding summand separately, i.e.
let us look at---no summation over $i$---

\begin{equation}\lae{9.43}
-(1-\frac\de2)\sum_k a^{ki} h^2_{kki} +\frac12 a^{kl} h_{kki} h_{lli} + \frac1{2
\widetilde H} \tilde F^{ii} (\sum_k h_{kki})^2,
\end{equation}

\cvm
\nd where we have divided the terms by $2$.

\cvm
Denote by $\sum'$ a sum where the index $i$ is omitted during the summation,
then, \re{9.43} can be expressed as

\begin{equation}\lae{9.44}
\begin{aligned}
-(1-\frac\de2)\sum_k a^{ki} h^2_{kki}&+h_{iii} \sum_k a^{ki} h_{kki}
+{\sum_{k<l}}' h_{kki} h_{lli} \\[\cma]
&+\frac1{2
\widetilde H} \tilde F^{ii} (\sum_k h_{kki})^2.
\end{aligned}
\end{equation}

\cvm
To replace $h_{iii}$ in the preceding expression we use the chain rule

\begin{equation}
f_i\equiv \tilde F_i= \tilde F^{kk} h_{kki}
\end{equation}

\cvm
\nd to derive

\begin{equation}\lae{9.46}
h_{iii}=\frac1{ \tilde F^{ii}}(f_i-{\sum_k}' \tilde F^{kk} h_{kki}).
\end{equation}

\cvm
Inserting \re{9.46} in \re{9.44} we obtain, after some simple algebraic
manipulations, cf. \ci[equ. (36) on p. 78]{pb} or \ci[equ. (17)]{pb2},

\begin{equation}\lae{9.47}
\begin{aligned}
&-(1-\frac\de2){\sum_k}' h^2_{kki}-{\sum_k}'\Big[\frac{ \tilde F^k_k}{ \tilde
F^i_i}-\frac{ \tilde F^i_i}{2 \widetilde  H}\Big(1-\frac{ \tilde F^k_k}{ \tilde
F^i_i}\Big)^2\Big]h^2_{kki}\\[\cma]
&- {\sum_{k<l}}'\Big[\frac{ \tilde F^k_k+ \tilde F^l_l}{ \tilde F^i_i}-1-\frac {
\tilde F^i_i}{ \widetilde H}\Big(1-\frac{ \tilde F^k_k}{ \tilde
F^i_i}\Big)\Big(1-\frac{ \tilde F^l_l}{ \tilde F^i_i}\Big)\Big]h_{kki}
h_{lli}\\[\cma]
&+{\sum_k}'\Big[\frac{ f_i}{ \tilde
F^i_i}+\frac{ f_i}{ \widetilde  H}\Big(1-\frac{ \tilde F^k_k}{ \tilde
F^i_i}\Big)\Big]h_{kki} + \frac{ \tilde F^i_i}{2 \widetilde H}\Big(\frac{f_i}{
\tilde F^i_i}\Big)^2.
\end{aligned}
\end{equation}

\cvm
Let us write \re{9.47} as the sum of three expressions $I_1+I_2+I_3$, where

\begin{equation}
I_1=-\frac\de2{\sum_k}'h^2_{kki} +{\sum_k}'\Big[\frac{ f_i}{ \tilde
F^i_i}+\frac{ f_i}{ \widetilde  H}\Big(1-\frac{ \tilde F^k_k}{ \tilde
F^i_i}\Big)\Big]h_{kki},
\end{equation}

\begin{equation}
I_2=\frac{ \tilde F^i_i}{2 \widetilde H}\Big(\frac{f_i}{ \tilde F^i_i}\Big)^2,
\end{equation}
\nd and
\begin{equation}
\begin{aligned}
I_3&= -(1-\de){\sum_k}'h^2_{kki}-{\sum_k}'\Big[\frac{ \tilde F^k_k}{ \tilde
F^i_i}-\frac{ \tilde F^i_i}{2 \widetilde  H}\Big(1-\frac{ \tilde F^k_k}{ \tilde
F^i_i}\Big)^2\Big]h^2_{kki}\\[\cma]
&- {\sum_{k<l}}'\Big[\frac{ \tilde F^k_k+ \tilde F^l_l}{ \tilde F^i_i}-1-\frac {
\tilde F^i_i}{ \widetilde H}\Big(1-\frac{ \tilde F^k_k}{ \tilde
F^i_i}\Big)\Big(1-\frac{ \tilde F^l_l}{ \tilde F^i_i}\Big)\Big]h_{kki}
h_{lli}.
\end{aligned}
\end{equation}

\cvm
In view of \re{9.23} we can estimate $I_1$ from above by

\begin{equation}\lae{9.51}
I_1\le c\msp[2]\de^{-1} f^{-1} \abs{f_i}^2.
\end{equation}

\cvm
$I_2$ is estimated by

\begin{equation}\lae{9.52}
I_2=\frac1{2 \widetilde H \tilde F^i_i}\abs{f_i}^2\le c f^{-1} \abs{f_i}^2,
\end{equation}

\cvm
\nd because of \re{1.17}.

\cvm
Finally, we claim that $I_3\le 0$ if we choose $\de=\frac14$. To verify this
assertion, let us multiply
$I_3$ by $2 \widetilde H \tilde F^i_i$ to obtain

\begin{equation}\lae{9.53}
\begin{aligned}
&-2(1-\de) \widetilde H \tilde F^i_i{\sum_k}'h^2_{kki}-{\sum_k}'[2 \widetilde H
\tilde F^k_k-( \tilde F^i_i- \tilde F^k_k)^2]h^2_{kki}\\[\cma]
&-{\sum_{k<l}}'[2 \widetilde H( \tilde F^k_k+ \tilde F^l_l)-2 \widetilde H \tilde
F^i_i-2 ( \tilde F^i_i- \tilde F^k_k)( \tilde F^i_i- \tilde F^l_l)] h_{kki} h_{lli}.
\end{aligned}
\end{equation}

\cvm
Now, we use \re{8.7} and replace \tit{any} $ \tilde F^j_j$, $1\le j\le n$, by

\begin{equation}
F^j_j+\e(n-1)(1+\e n) Hg^j_j\equiv F^j_j +\e\ga_\e H.
\end{equation}

\cvm
The expression in \re{9.53} is then equal to the sum of two terms $I_4+I_5$,
where

\begin{equation}\lae{9.55}
\begin{aligned}
&I_4=
-2(1-\de) \widetilde H  F^i_i{\sum_k}'h^2_{kki}-{\sum_k}'[2
\widetilde H
 F^k_k-(  F^i_i-  F^k_k)^2]h^2_{kki}\\[\cma]
&-{\sum_{k<l}}'[2 \widetilde H(  F^k_k+  F^l_l)-2 \widetilde H 
F^i_i-2 (  F^i_i-  F^k_k)(  F^i_i-  F^l_l)] h_{kki} h_{lli},
\end{aligned}
\end{equation}

\cvm
\nd and

\begin{multline}
I_5= -2 (1-\de) \widetilde H\e\ga_\e H {\sum_k}'h^2_{kki} - 2\widetilde H \e
\ga_\e H {\sum_k}' h^2_{kki}\\[\cma]
 - 2 \widetilde H \e \ga_\e H{\sum_{k<l}}' h_{kki}h_{lli}. 
\end{multline}

From the the binomial formula 

\begin{equation}
\Big({\sum_k}'h_{kki}\Big)^2={\sum_k}'h^2_{kki}+ 2{\sum_{k<l}}' h_{kki}
h_{lli}
\end{equation}

\cvm
\nd we infer that $I_5\le 0$, while $I_4$ is non-positive provided we choose
$\de=\frac14$ and assume

\begin{equation}\lae{9.58}
\max(- \tilde \ka_1,0)\le \frac1{2(n-1)} \widetilde  H,
\end{equation}

\cvm
\nd cf. \ci[pp. 81--85]{pb} or \ci[pp. 23--29]{pb2}.

 But the condition \re{9.58} is certainly satisfied in
view of
\re{9.28}.

Combining \re{9.30}, \re{9.42}, \re{9.51}, and \re{9.52} gives \re{9.27}, and
thus,  the Lemma is proved.
\ep

\cvb
From \rl{9.2} and \rl{9.5} we conclude

\cvb
\bc\lac{9.6}
Let $M=\graph u$ be an admissible solution of equation \re{9.1} in $\Om$. Then,
the estimate

\begin{multline}\lae{9.59}
\tilde F^{ij,kl}h_{ij;r}h_{kl;s} g^{rs}H^{-1} + \tilde F^{ij} (\log H)_i (\log
H)_j\\
\le c\msp H^{-1} f^{-1}
\norm{Df}^2 + c \msp H^{-1} \norm{D\log H}^2 + c H^{-1}
\end{multline}

\cvm
\nd is valid in every point $p\in M$, where \re{9.28} is satisfied. The constant
$c$ depends on $\Om$, $\nnorm{Df}, \nnorm{D^2f}$,  the
constant $c_1$ in \re{8.3}, and on known estimates of the $C^0$ and $C^1$-
norm of $u$.
\ec

\cvb
As we already mentioned we have to assume the existence of a strictly convex
function $\chi\in C^2(\bar\Om)$, i.e. $\chi$ satisfies

\begin{equation}
\chi_{\al\bet}\ge c_0\msp  \bar g_{\al\bet}
\end{equation}

\cvm
\nd with a positive constant $c_0$.

We observe that then the restriction $\chi=\fv\chi M$ of $\chi$ to an admissible
solution $M\su \bar\Om$ of \re{9.1} satisfies the elliptic inequality

\begin{equation}\lae{9.61}
\begin{aligned}
- \tilde F^{ij}\chi _{{\vphantom{}}_{ij}}&=-2 \tilde F\chi _{{\vphantom{}}_\al}
\n^\al - \tilde F^{ij}\chi _{{\vphantom{}}_{\al\bet}} x^\al_i x^\bet_j\\[\cma]
&\le -2 \tilde F\chi _{{\vphantom{}}_\al}
\n^\al -c_0 \tilde F^{ij}g_{ij},
\end{aligned}
\end{equation}

\cvm
\nd where we used the homogeneity of $ \tilde F$.

\cvb
We can now prove uniform $C^2$- estimates.

\cvb

\bt
Let $M=\graph u$ be an admissible solution of equation \re{9.1} in $\Om$,
where $f$ satisfies the estimates \re{8.3}, \re{8.4} and \re{8.5}. Then, the
principal curvatures of $M$ are uniformly bounded.
\et

\bp
Let $\chi$ be the strictly convex function and $\m$ a large positive constant.
We shall prove that $w=\log H +\m\chi$ is uniformly bounded from above.

\cvm
Let $x_0\in M$ be such that

\begin{equation}
w(x_0)=\sup_M w,
\end{equation}

\cvm
\nd and choose in $x_0$ a local coordinate system satisfying \re{8.14} and
\re{8.15}. Applying the maximum principle, we conclude from \re{9.3} and
\re{9.61}

\begin{equation}\lae{9.63}
\begin{aligned}
0&\le - \tilde F^{kl} h_{kr} h^r_l+c \tilde F^{ij} g_{ij} +c \m f -\m c_0 \tilde
F^{ij} g_{ij}\\[\cma]
&\hp{\le \; }+ c(1+f+\nnorm{Df}+\nnorm{D^2f})(1+H+\norm{D\log H})\\[\cma]
&\hp{\le \; }+\tilde F^{ij,kl}h_{ij;r}h_{kl;s} g^{rs}H^{-1} + \tilde F^{ij} (\log
H)_i (\log H)_j,
\end{aligned}
\end{equation}

\cvm
\nd where we also assumed $H$ to be larger than $1$.

\cvb
We now consider two cases.

\cvb
\tit{Case $1$}: Suppose that

\begin{equation}
\abs{\ka_1}\ge \e_1 H\equiv \frac1{2(n-1)} H.
\end{equation}

\cvm
Then, we infer from \rl{8.3} and \re{9.24}

\begin{equation}
- \tilde F^{kl} h_{kr} h^r_l\le -\frac{n-1}n H \ka^2_1\le - \frac{n-1}n \e^2_1
H^3 \equiv -\e_2 H^3.
\end{equation}

\cvm
Moreover, the concavity of $\log \tilde F$ implies

\begin{equation}
\begin{aligned}
\tilde F^{ij,kl}h_{ij;r}h_{kl;s} g^{rs}H^{-1}&\le \tilde F^{-1}g^{ij} \tilde F^{kl}
h_{kl;i}
\tilde F^{rs} h_{rs;j} H^{-1}\\[\cma]
&= f^{-1} \norm{Df}^2 H^{-1}\\[\cma]
&\le c f^{-1}\nnorm{Df}^2 \abs A^2H^{-1}\\[\cma]
&\le c f^{-1} \nnorm{Df}^2 H.
\end{aligned}
\end{equation}

\cvm
Furthermore, $Dw(x_0)=0$, or,

\begin{equation}\lae{9.67}
(\log H)_i=-\m \chi _{{\vphantom{}}_i}.
\end{equation}

\cvm
Inserting the last three relations in \re{9.63} we obtain

\begin{equation}
0\le -\e_2 H^3 +c(1+H+\m) + c\m^2 H,
\end{equation}

\cvm
\nd where, now, $c$ depends on $f$ and its derivatives in the ambient space.

Hence, $H$, and therefore $w$, are a priori bounded in $x_0$.

\cvb
\tit{Case $2$}: Suppose that

\begin{equation}
\abs{\ka_1}< \e_1 H.
\end{equation}

\cvm
Then, \rc{9.6} is applicable, and we infer from \re{9.63} and \re{9.67}

\begin{equation}
0\le c (1+H+\m+\m^2 H^{-1}) + (c-\m c_0) \tilde F^{ij}g_{ij}.
\end{equation}

\cvm
Choosing now $\m$ sufficiently large we obtain an a priori bound for $H(x_0)$,
since

\begin{equation}
\tilde F^{ij}g_{ij}\ge (n-1) H.
\end{equation}

\cvm
Thus, $w$, or equivalently $H$, are uniformly bounded.
\ep

\cvb
\section{Existence of a solution}\las{10}

\cvb
We can now demonstrate the final step in the proof of \rt{0.2}. Let $M_\e=\graph
u_\e$ be the stationary approximations. In the preceding sections we have
proved uniform estimates for $u_\e$ up to the order two. Since, by assumption,
$f$ is strictly positive, the principal curvatures of
$M_\e$ stay in a compact subset of the cone $\C_2$ for small $\e$, cf. \rr{3.1}, 
and therefore, the operator
$
\tilde F$ is uniformly elliptic for those $\e$. Taking the square root on both sides
of equation
\re{9.1} without changing the notation, we also know that $ \tilde F$ is concave.

Hence the $C^{2,\al}$- estimates of Evans and Krylov are applicable, cf. \ci{ce}
and \ci{nk}, and we deduce

\begin{equation}
\abs{u_\e}_{2,\al,\so}\le \const
\end{equation}

\cvm
\nd uniformly in $\e$. If $\e$ tends to zero, a subsequence converges to a
solution $u\in C^{2,\al}(\so)$ of our problem. From the Schauder estimates we
further conclude $u\in C^{4,\al}(\so)$.

\cvb

\end{document}